\documentclass[AMA,Times1COL]{WileyNJDv5}
\usepackage{subcaption}

\articletype{Research Article}%

\received{}
\revised{}
\accepted{}
\journal{International Journal for Numerical Methods in Engineering}
\volume{}
\copyyear{}
\startpage{}

\begin{document}

\title{Super-Resolution of Elliptic PDE Solutions Using Least Squares Support Vector Regression}

\author[1]{Maryam Babaei}
\author[3]{Péter Rucz}
\author[1]{Manfred Kaltenbacher}
\author[1]{Stefan Schoder}

\authormark{M. Babaei, \textsc{et al.}}
\titlemark{Super-Resolution of Elliptic PDE Solutions Using LS-SVR}

\address[1]{\orgdiv{Institute of Fundamentals and Theory in Electrical Engineering}, \orgname{Graz University of Technology}, \orgaddress{\state{Graz}, \country{Austria}}}


\address[3]{\orgdiv{Department of Networked Systems and Services, Faculty of Electrical Engineering and Informatics}, \orgname{Budapest University of Technology and Economics}, \orgaddress{\state{Budapest}, \country{Hungary}}}

\corres{Corresponding author Stefan Schoder  \email{ stefan.schoder@tugraz.at}}

\abstract[Abstract]{%
A hybrid computational approach that integrates the finite element method (FEM) with least squares support vector regression (LSSVR) is introduced to solve partial differential equations. The method combines FEM's ability to provide the nodal solutions and LSSVR with higher-order Legendre polynomial kernels to deliver a closed-form analytical solution for interpolation between the nodes. The hybrid approach implements element-wise enhancement (super-resolution) of a given numerical solution, resulting in high resolution accuracy, while maintaining consistency with FEM nodal values at element boundaries. It can adapt any low-order FEM code to obtain high-order resolution by leveraging localized kernel refinement and parallel computation without additional implementation overhead. Therefore, effective inference/post-processing of the obtained super-resolved solution is possible. Evaluation results show that the hybrid FEM-LSSVR approach can achieve significantly higher accuracy compared to the base FEM solution. Comparable accuracy is a achieved when comparing the hybrid solution with a standalone FEM result with the same polynomial basis function order. The convergence studies were conducted for four elliptic boundary value problems to demonstrate the method's ability, accuracy, and reliability. Finally, the algorithm can be directly used as a plug-and-play method for super-resolving low-order numerical solvers and for super-resolution of expensive/under-resolved experimental data.}


\maketitle

\renewcommand\thefootnote{}
\footnotetext{\textbf{Abbreviations:} FEM, finite element method; LSSVR, least squares support vector regression; PDE, partial differential equation; MAE, mean absolute error; MSE, mean squared error.}

\renewcommand\thefootnote{\fnsymbol{footnote}}
\setcounter{footnote}{1}

\section{Introduction}
Partial differential equations are used to describe structural mechanics phenomena across various fields of science and engineering \cite{kaltenbacher2015numerical, bies2023engineering,rufo2022sound}, encompassing acoustics, fluid dynamics, electromagnetics, and coupled fields. For instance, at low excitation frequency, the room acoustics behavior is mainly governed by the modal wave field, being the solution of the wave equation \cite{pierce2019acoustics}. Engineers are solving these equations when predicting the acoustic behavior of products, designing sound systems, and developing noise reduction technologies. For numerical methods, achieving high-accuracy solutions is computationally demanding.
Additionally, key challenges arise in parameter estimation and inverse problems, where unknown material properties or source terms must be determined from limited measurements \cite{bryan2020impulse, muller2024free}. Working with sparse or noisy data makes it harder to uncover the underlying physics, especially when traditional interpolation methods fail to capture the underlying physics between measurement points \cite{mezza2024data}. Machine learning can support this process by combining data-driven approaches with physical constraints \cite{raissi2019physics}.


The combination of machine learning with numerical methods for solving PDEs has emerged recently. Physics-informed neural networks (PINNs) have demonstrated solving forward and inverse problems by embedding physical laws directly into neural network loss functions for training \cite{raissi2019physics,karniadakis2021physics,schoder2025physics,mirzabeigi2025physics}. Recent advancements have improved these methods by developing specialized network architectures and enhanced training strategies to address convergence issues and improve solution accuracy \cite{lu2021deepxde, cuomo2022scientific, mohammadi2025new}. In loss-based PDE and boundary condition enforcement, the different loss components compete during training and require careful hyperparameter tuning to balance the weighting of different loss components \cite{cuomo2022scientific,sharma2023stiff}. These loss-based formulations offer no guarantee of finding valid PDE solutions and may suffer from local minima issues inherent in neural network optimization. An alternative approach would be for the optimizer to take care of the PDE constraint; however, such optimizers are typically challenging to develop and require significant solver expertise \cite{mehrkanoon2015learning}. In contrast, the LSSVR approach is an automated method that treats the PDE as an equality constraint in the optimization, directly utilizing a linear PDE as a constraint in the optimizer routine. This structured method eliminates the high constraint-optimizer implementation costs and the need for an additional hyperparameter that balances a PDE loss against a data loss of PINNs. Global optimality and PDE satisfaction is achieved by transforming the problem into a convex quadratic programming form. This is fundamentally different from loss-based approaches in combination with iterative optimization methods, which often face local minima issues and require careful hyperparameter tuning to balance competing objectives. 
LSSVR supports various kernel functions, including radial basis functions, Gaussian kernels, polynomial kernels, and different orthogonal polynomial families. The least squares support vector machine method \cite{mehrkanoon2012approximate,mehrkanoon2015learning} and its support of higher-order (polynomial) kernel functions offers a perspective to reduce this interpolation error \cite{padierna2018novel, moghaddam2016new, babaei2025physics}. For instance, Legendre polynomials form an orthogonal basis on $[-1,1]$ with respect to the inner product $\langle f,g \rangle = \int_{-1}^{1} f(x)g(x)\,\mathrm{d}x$. When using this on a real physical domain, we need to apply a simple linear scaling to map each element to the standard interval $[-1, 1]$. Using these polynomials as kernel functions, theoretically leads to spectral convergence rates for smooth functions \cite{castillo2024first}. In spectral element methods, Legendre polynomials serve as basis functions that allow accurate solutions with fewer degrees of freedom compared to standard FEM. When using Legendre polynomials as kernel functions in LSSVR, they enable the solution to adapt to complex behaviors between interpolation points, thereby improving accuracy. This property leads to a fundamental advantage of this approach, which lies in its closed-form solution representation and allows for efficient predictions. 

This work investigates three hypotheses regarding the combination of FEM and LSSVR into a hybrid FEM-LSSVR method for solving partial differential equations. 
The first hypothesis posits that a hybrid FEM-LSSVR method can be established for solving elliptic boundary value problems. It will utilize standard libraries for LSSVR Kernel refinement on an element, allowing the computation of derivatives of the solution via automatic differentiation. This localized method allows us to compute derivatives for post-processing and error estimation based on the PDE residual.

The second hypothesis is that the hybrid FEM-LSSVR approach can significantly improve solution accuracy between nodal values provided from any source, including low-order polynomial basis function finite element method results, experimental measurement data, or output from any kind of boundary value solvers. This flexibility to incorporate external information allows us to compute a super-resolved field by increasing the effective resolution of the solution between the nodes. It is also hypothesized that this approach can overcome the degraded convergence that low-order polynomial basis functions used in FEM solvers exhibit by combining it with an automatic method to increase the order of the interpolation kernel, as implemented in the LSSVR method.

Thirdly and finally, the smoothness hypothesis is that the FEM-LSSVR method, by using high-order Legendre polynomial kernels, yields a solution field with increased local smoothness and enhanced spatial derivatives. This property is crucial for derivative-intensive post-processing applications (e.g., computational aeroacoustics \cite{Schoder2022b,schoder2020aeroacoustic2} or post-processing fluid flows \cite{schoder2020postprocessing,krais2021flexi}). 
%
%




The rest of the article is organized as follows: Section~\ref{sec:methodology} introduces the theoretical background of the FEM and presents the LSSVR formulation with Legendre kernels for solving elliptic equations, including both Poisson and Helmholtz equations, and detailing the hybrid FEM–LSSVR algorithm. Section~\ref{sec:results} displays the numerical results for different test problems, including convergence studies and performance comparisons for one- and two-dimensional Poisson and Helmholtz equations with different boundary conditions. Finally, Section~\ref{sec:Conclusion} concludes the work and suggests directions for future research.
\section{Methodology}
\label{sec:methodology}
In this section, we provide a brief introduction to the FEM and present the theory of the LSSVR and the hybrid FEM-LSSVR approach. We explain the application of LSSVR with Legendre polynomial kernels and show how the system matrices are assembled by the algorithm.

\subsection{Finite Element Method}\label{sec:fem}

The FEM aims to find the best approximate solutions to boundary value problems for PDEs under a suitable norm. Thereby, a strong formulation of a partial differential equation is transformed into a weak form through the process of partial integration, and then the Galerkin variational method is applied to obtain an algebraic system of equations \cite{kaltenbacher2015numerical}.

In this article, we consider the Poisson equation
\begin{equation}\label{eq:poisson}
-\nabla^2u(x,y)=f(x,y), \quad (x,y)\in\Omega
\end{equation}
and the Helmholtz equation
\begin{equation}\label{eq:helmholtz}
\nabla^2u+\xi ^2u=f(x,y), \quad (x,y)\in\Omega.
\end{equation}
where $\xi$ represents the wavenumber.
The boundary conditions for these elliptic equations include Dirichlet boundary conditions
\begin{equation}
u(x,y)=g(x,y), \quad (x,y)\in\Gamma_D \subset \partial\Omega
\end{equation}
and Neumann boundary conditions
\begin{equation}
\frac{\partial u}{\partial n}(x,y)=h(x,y), \quad (x,y)\in\Gamma_N \subset \partial\Omega,
\end{equation}
where the domain boundary satisfies
\begin{equation}
\Gamma_D \cup \Gamma_N = \partial\Omega
\end{equation}
and $\partial \Omega$ denotes the boundary of the domain $\Omega$.
For the Poisson equation, the weak formulation becomes
\begin{equation}
\int_\Omega \nabla u \cdot \nabla v \, \mathrm{d}\Omega = \int_\Omega fv \, \mathrm{d}\Omega,
\end{equation}
while for the Helmholtz equation, it takes the form
\begin{equation}
-\int_\Omega \nabla u \cdot \nabla v \, \mathrm{d}\Omega + \xi^2\int_\Omega uv \, \mathrm{d}\Omega = \int_\Omega fv \, \mathrm{d}\Omega,
\end{equation}
with $v$ being a scalar test function.
In FEM, the domain $\Omega$ is partitioned into a mesh of finite elements, and the solution $u$ is approximated as
\begin{equation}
u(x,y) \approx u_h(x,y) = \sum_{j=1}^{N} u_j \phi_j(x,y),
\end{equation}
where $\phi_j$ are piecewise polynomial basis functions of a chosen order and $u_j$ are unknown coefficients.
Substituting this approximation into the weak form and choosing test functions $v = \phi_i$ (Galerkin method), the linear system becomes
\begin{equation}
\boldsymbol{\mathcal{K}} \boldsymbol{u} + \boldsymbol{\mathcal{M}} \boldsymbol{u} = \boldsymbol{f}.
\end{equation}
The system matrices and load vector are defined as
\begin{align*}
\mathcal{K}_{ij} &= \int_\Omega \nabla \phi_j \cdot \nabla \phi_i \, \mathrm{d}\Omega \quad \text{(stiffness matrix)} \\
\mathcal{M}_{ij} &= \xi^2\int_\Omega \phi_j \phi_i \, \mathrm{d}\Omega \quad \text{(mass matrix, Helmholtz only)} \\
f_i &= \int_\Omega f \phi_i \, \mathrm{d}\Omega \quad \text{(load vector)}.
\end{align*}
For the FEM implementation and evaluation, Scikit-FEM 11.0.0 \cite{skfem2020} in Python 3.11.13 is used.
\subsection{LSSVR for Solving Elliptic Equations}

The LSSVR approach is developed to solve elliptic partial differential equations with different boundary conditions. In this section, we apply LSSVR to equations \eqref{eq:poisson} and \eqref{eq:helmholtz}. The goal is to obtain an approximate solution $\hat{u}(x,y)$ of an elliptic equation using the LSSVR framework. The method will find the solution by minimizing errors in the differential equation while avoiding overfitting through a regularization term.
The approximate solution is assumed to have the form of
\begin{equation*}
\hat{u}(x,y) = \mathbf{w}^\mathrm{T}\boldsymbol{\phi}(x,y) + b,
\end{equation*}
where $\mathbf{w} \in \mathbb{R}^h$ is a weight vector, $b \in \mathbb{R}$ is a bias term, and $\boldsymbol{\phi}: \mathbb{R}^2 \to \mathbb{R}^h$ represents the feature mapping induced by the Legendre polynomial kernel into an $h$-dimensional feature space. The relationship between the feature mapping $\boldsymbol{\phi}$ and the kernel function follows Mercer's theorem through \cite{scholkopf2002learning}
\begin{equation*}
K((x,y),(x',y')) = \boldsymbol{\phi}(x,y)^T\boldsymbol{\phi}(x', y'),
\end{equation*}
which allows the implicit use of the feature space through the kernel without explicit computation of $\boldsymbol{\phi}$.

Following the LSSVR framework for differential equations \cite{mehrkanoon2012approximate}, the optimization problem is formulated to enforce the differential equation at interior collocation points while satisfying boundary conditions. For the general elliptic equation case, let $\mathcal{L}$ represent the differential operator, where $\mathcal{L} = -\nabla^2$ for the Poisson equation and $\mathcal{L} = \nabla^2 + \xi^2$ for the Helmholtz equation. The constraints ensure that the approximate solution satisfies the elliptic equation at $N$ selected interior points $(x_i,y_i)$ within the domain, while boundary constraints enforce the conditions on $\partial\Omega$
\begin{equation}
\begin{aligned}
\min_{\mathbf{w},b,\mathbf{e}} &\quad \frac{1}{2}\mathbf{w}^T\mathbf{w}+\frac{\gamma}{2}\sum_{i=1}^{N}e_i^2 \\ 
\text{s.t.} &\quad \mathcal{L}(\mathbf{w}^T\boldsymbol{\phi}(x_i,y_i)+b) = f(x_i,y_i)+e_i, \quad i=1,\ldots,N, \\ 
&\quad \mathbf{w}^T\boldsymbol{\phi}(x,y)+b = g(x,y), \quad (x,y)\in\Gamma_D, \\
&\quad \frac{\partial}{\partial n}(\mathbf{w}^T\boldsymbol{\phi}(x,y)+b) = h(x,y), \quad (x,y)\in\Gamma_N.
\end{aligned}
\end{equation}
Here, $\gamma > 0$ is a regularization parameter, and $\mathbf{e} = [e_1, \ldots, e_N]^T$ is the vector of the residual of the PDE at the $N$ interior collocation points $(x_i,y_i)$. The error variables $e_i$ in the LSSVR formulation are essential for addressing the fundamental mismatch between the infinite-dimensional space where PDE solutions naturally exist and the finite-dimensional kernel-induced approximation space \cite{Evgeniou2000,suykens2002least}. Without these slack variables, the optimization problem becomes infeasible when the finite-dimensional kernel space cannot exactly represent the PDE solution, particularly for high-frequency oscillatory problems where polynomial kernels of limited degree fail to capture rapid variations \cite{Tikhonov1977,Engl1996}. The regularization parameter $\gamma$ is controlling the penalty on the residual $e_i$. It balances between constraint satisfaction and solution smoothness, ensures numerical stability and prevents overfitting through Tikhonov regularization principles \cite{Scholkopf2002,Wahba1990}.

The Lagrangian for this constrained optimization problem with mixed boundary conditions is formulated by
\begin{align}
\mathbf{L}(\mathbf{w}, b, \mathbf{e}, \boldsymbol{\alpha}, \boldsymbol{\beta}_D, \boldsymbol{\beta}_N) &= \frac{1}{2}\mathbf{w}^T\mathbf{w}+\frac{\gamma}{2}\sum_{i=1}^{N}e_i^2 \nonumber \\
&\quad - \sum_{i=1}^{N} \alpha_i \left[\mathcal{L}(\mathbf{w}^T\boldsymbol{\phi}(x_i,y_i)+b) - f(x_i,y_i) - e_i\right] \nonumber \\
&\quad - \sum_{j=1}^{N_D} \beta_{D,j} \left[\mathbf{w}^T\boldsymbol{\phi}(x_j^D,y_j^D)+b - g(x_j^D,y_j^D)\right] \nonumber \\
&\quad - \sum_{k=1}^{N_N} \beta_{N,k} \left[\frac{\partial}{\partial n}(\mathbf{w}^T\boldsymbol{\phi}(x_k^N,y_k^N)+b) - h(x_k^N,y_k^N)\right],
\end{align}
where $\boldsymbol{\alpha} = [\alpha_1, \ldots, \alpha_N]^T$, $\boldsymbol{\beta}_D = [\beta_{D,1}, \ldots, \beta_{D,N_D}]^T$, and $\boldsymbol{\beta}_N = [\beta_{N,1}, \ldots, \beta_{N,N_N}]^T$ are the vectors of Lagrange multipliers corresponding to the PDE constraints, Dirichlet boundary conditions, and Neumann boundary conditions, respectively.
The kernel function is defined by using Legendre polynomials as
\begin{equation}
K((x,y),(x',y')) = \sum_{i=0}^{M_x-1}\sum_{j=0}^{M_y-1} P_i(x)P_j(y)P_i(x')P_j(y'),
\end{equation}
where $P_i$ and $P_j$ are Legendre polynomials of orders $i$ and $j$ respectively, and $M_x$, $M_y$ define the polynomial orders in the $x$ and $y$ directions. The kernel choice can provide a natural compatibility with FEM basis functions of the same polynomial order. While the higher order terms can improve the method's interpolation capabilities. The kernel matrix elements are computed as
\begin{equation*}
K_{ij}=K((x_i,y_i),(x_j,y_j)).
\end{equation*}
Applying the Karush-Kuhn-Tucker (KKT) optimality conditions to the Lagrangian leads to
\begin{align}
\frac{\partial \mathbf{L}}{\partial \mathbf{w}} &= \mathbf{w} + \sum_{i=1}^{N} \alpha_i \mathcal{L}\boldsymbol{\phi}(x_i,y_i) + \sum_{j=1}^{N_D} \beta_{D,j} \boldsymbol{\phi}(x_j^D,y_j^D) + \sum_{k=1}^{N_N} \beta_{N,k} \frac{\partial}{\partial n}\boldsymbol{\phi}(x_k^N,y_k^N) = 0 \\
\frac{\partial \mathbf{L}}{\partial b} &= \sum_{i=1}^{N} \alpha_i + \sum_{j=1}^{N_D} \beta_{D,j} + \sum_{k=1}^{N_N} \beta_{N,k} = 0 \\
\frac{\partial \mathbf{L}}{\partial e_i} &= \gamma e_i + \alpha_i = 0, \quad i = 1,\ldots,N.
\end{align}
Following the approach in \cite{mehrkanoon2015learning}, the elimination of primal variables and application of Mercer's theorem yields the dual formulation as a linear system
\begin{equation}
\begin{bmatrix} 
\mathbf{K} + \gamma^{-1}\mathbf{I} & \mathbf{S} & \mathbf{b} \\
\mathbf{S}^T & \boldsymbol{B} & \mathbf{1}_{M} \\
\mathbf{b}^T & \mathbf{1}_{M}^T & 0
\end{bmatrix} 
\begin{bmatrix} 
\boldsymbol{\alpha} \\
\boldsymbol{\beta} \\
b
\end{bmatrix} 
= 
\begin{bmatrix} 
\mathbf{f} \\
\mathbf{v} \\
0
\end{bmatrix}.
\end{equation}
The matrix $\mathbf{K} \in \mathbb{R}^{N \times N}$ results from applying the differential operator $\mathcal{L}$ to the kernel function at interior collocation points. For the Poisson equation in two dimensions, this becomes
\begin{equation*}
K_{ij} = -\nabla^2 K((x_i,y_i),(x_j,y_j)) = -\left[\frac{\partial^2 K}{\partial x_j^2} + \frac{\partial^2 K}{\partial y_j^2}\right]_{(x_i,y_i),(x_j,y_j)},
\end{equation*}
while for the Helmholtz equation, it is 
\begin{equation*}
K_{ij} = (\nabla^2 + \xi^2) K((x_i,y_i),(x_j,y_j)).
\end{equation*}
The coupling matrix $\mathbf{S} \in \mathbb{R}^{N \times M}$ is formed by concatenating the Dirichlet and Neumann coupling matrices as $\mathbf{S} = [\mathbf{S}_D, \mathbf{S}_N]$, where $\mathbf{S}_D \in \mathbb{R}^{N \times N_D}$ represents the coupling between interior points and Dirichlet boundary points through kernel evaluations, and $\mathbf{S}_N \in \mathbb{R}^{N \times N_N}$ represents the coupling between interior points and Neumann boundary points through normal derivative evaluations of the kernel function. The specific entries are computed as
\begin{align*}
[\mathbf{S}_D]_{ij} &= K((x_i, y_i), (x_j^D, y_j^D)), \\
[\mathbf{S}_N]_{ij} &= \frac{\partial K((x_i, y_i), (x_j^N, y_j^N))}{\partial n_j},
\end{align*}
where $(x_i, y_i)$ for $i = 1, \ldots, N$ are interior collocation points, $(x_j^D, y_j^D)$ for $j = 1, \ldots, N_D$ are Dirichlet boundary points, $(x_j^N, y_j^N)$ for $j = 1, \ldots, N_N$ are Neumann boundary points, and $\frac{\partial}{\partial n_j}$ denotes the outward normal derivative at the $j$-th Neumann boundary point.
The boundary kernel matrix $\boldsymbol{B} \in \mathbb{R}^{M \times M}$ has the block structure
\begin{equation*}
\boldsymbol{B} = \begin{bmatrix} \boldsymbol{B}_D & \boldsymbol{B}_{DN} \\ \boldsymbol{B}_{DN}^T & \boldsymbol{B}_N \end{bmatrix}
\end{equation*}
where $\boldsymbol{B}_D \in \mathbb{R}^{N_D \times N_D}$ contains kernel evaluations between Dirichlet boundary points, $\boldsymbol{B}_N \in \mathbb{R}^{N_N \times N_N}$ contains normal derivative evaluations between Neumann boundary points, and $\boldsymbol{B}_{DN} \in \mathbb{R}^{N_D \times N_N}$ contains cross-evaluations between Dirichlet and Neumann boundary points. The entries are calculated as
\begin{align*}
[\boldsymbol{B}_D]_{ij} &= K((x_i^D, y_i^D), (x_j^D, y_j^D)), \\
[\boldsymbol{B}_N]_{ij} &= \frac{\partial^2 K((x_i^N, y_i^N), (x_j^N, y_j^N))}{\partial n_i \partial n_j}, \\
[\boldsymbol{B}_{DN}]_{ij} &= \frac{\partial K((x_i^D, y_i^D), (x_j^N, y_j^N))}{\partial n_j}.
\end{align*}
The Lagrangian multiplier vector is structured as $\boldsymbol{\beta} = [\boldsymbol{\beta}_D^T, \boldsymbol{\beta}_N^T]^T \in \mathbb{R}^M$ where $\boldsymbol{\beta}_D \in \mathbb{R}^{N_D}$ corresponds to Dirichlet constraints and $\boldsymbol{\beta}_N \in \mathbb{R}^{N_N}$ corresponds to Neumann constraints. Similarly, the boundary condition vector is $\mathbf{v} = [\mathbf{g}^T, \mathbf{h}^T]^T \in \mathbb{R}^M$ containing both Dirichlet values $\mathbf{g} \in \mathbb{R}^{N_D}$ and Neumann values $\mathbf{h} \in \mathbb{R}^{N_N}$. The total number of boundary points is $M = N_D + N_N$, and $\mathbf{b} = [1, \ldots, 1]^T \in \mathbb{R}^N$ is the vector of ones.

After solving this linear system, the solution at any point $(x,y)$ in the domain is evaluated using
\begin{align}
u(x,y) &= \sum_{i=1}^{N} \alpha_i K((x,y),(x_i,y_i)) + \sum_{k=1}^{N_D} \beta_{D,k} K((x,y),(x_k^D,y_k^D)) \\
&\quad + \sum_{\ell=1}^{N_N} \beta_{N,\ell} K((x,y),(x_\ell^N,y_\ell^N)) + b,
\end{align}
where $(x_i,y_i)$ are the interior collocation points, $(x_k^D,y_k^D)$ are the Dirichlet boundary points, $(x_\ell^N,y_\ell^N)$ are the Neumann boundary points, and $\boldsymbol{\beta} = [\boldsymbol{\beta}_D^T, \boldsymbol{\beta}_N^T]^T$ contains the Lagrangian multipliers for all boundary conditions. This dual representation provides the closed-form solution that can be evaluated at any point in the domain without requiring additional interpolation procedures.
\subsection{Hybrid FEM-LSSVR Method}
The proposed Hybrid FEM-LSSVR method combines the FEM with least squares support vector regression. Based on the solution of the FEM method, the LSSVR with Legendre polynomial kernels aims to improve the interpolation accuracy. The FEM handles the boundary conditions and the resolution of the geometry. The hybrid methodology transforms an (initial) FEM solution into an enhanced global approximation through three main steps: First, we use standard FEM to get solutions at the nodes. Second, we improve the solution inside each element using LSSVR with Legendre polynomial kernels while keeping the FEM nodal values fixed. Third, we build the global solution from these improved elements. The FEM handles the boundary conditions and domain geometry, while LSSVR improves the accuracy between nodes.

For the first step of the algorithm, we apply the standard finite element method to discretize the domain $\Omega$ and obtain a preliminary solution at the nodal points. The FEM solution is given by
\begin{equation}
u^{\mathrm{fem}}(\mathbf{x}) = \sum_{i=1}^{n_{\mathrm{fem}}} u_i^{\mathrm{fem}} N_i(\mathbf{x}) ,
\end{equation}
where $\{\mathbf{x}_i^{\mathrm{fem}}\}_{i=1}^{n_{\mathrm{fem}}}$ are the node positions, $\{u_i^{\text{fem}}\}_{i=1}^{n_{\text{fem}}}$ are the corresponding solution values, and $N_i(\mathbf{x})$ are the finite element shape functions.

Step 2: Each finite element in the domain decomposition stage becomes a subdomain for localized LSSVR enhancement. We split the domain $\Omega$ into $E_{\mathrm{fem}}$ elements $\{\Omega_e\}_{e=1}^{E_{\mathrm{fem}}}$, so that $\Omega = \bigcup_{e=1}^{E_{\mathrm{fem}}} \Omega_e$.  For each element $\Omega_e$, we apply the LSSVR method with Legendre polynomial kernels. This sets up an optimization problem that improves the solution interpolation while keeping the FEM nodal values fixed along the element boundaries.
The optimization problem for each element is then defined as
\begin{equation}
\begin{aligned}
\min_{\mathbf{w}_e,b_e,\mathbf{e}_e} &\quad \frac{1}{2}\mathbf{w}_e^T\mathbf{w}_e + \frac{\gamma}{2}\sum_{j=1}^{M_e} e_{e,j}^2 \\
\text{s.t.} &\quad \mathcal{L}[\mathbf{w}_e^T\boldsymbol{\phi}(\mathbf{x}_j^{(e)}) + b_e] = f(\mathbf{x}_j^{(e)}) + e_{e,j}, \quad j=1,\ldots,M_e \\
&\quad \mathbf{w}_e^T\boldsymbol{\phi}(\mathbf{x}_k) + b_e = u_k^{\mathrm{fem}}, \quad \mathbf{x}_k \in \partial\Omega_e \cap \{\mathbf{x}_i^{\mathrm{fem}}\}_{i=1}^{n_{\mathrm{fem}}}
\end{aligned}
\end{equation}
where $\{\mathbf{x}_j^{(e)}\}_{j=1}^{M_e}$ are the training points inside each element $\Omega_e$, and $u_k^{\mathrm{fem}}$ are the FEM nodal values on the element boundaries $\partial \Omega_e$.
As shown in the previous section, by solving the LSSVR optimization problem we obtain a dual formulation on a finite element level
\begin{equation}
\begin{bmatrix}
\mathbf{K}_e + \gamma^{-1}\mathbf{I} & \mathbf{S}_{e} & \mathbf{b}_e \\ 
\mathbf{S}_{e}^T & \boldsymbol{B}_{e} & \mathbf{1}_{|\mathcal{I}_e|} \\
\mathbf{b}_e^T & \mathbf{1}_{|\mathcal{I}_e|}^T & 0
\end{bmatrix}
\begin{bmatrix}
\boldsymbol{\alpha}_e \\
\boldsymbol{\beta}_e \\
b_e
\end{bmatrix}
=
\begin{bmatrix}
\mathbf{f}_e \\
\mathbf{u}_e^{\text{fem}} \\
0
\end{bmatrix},
\label{eq:dual_system}
\end{equation}
where the matrix dimensions and computational entries are:
\begin{align}
\mathbf{K}_e &\in \mathbb{R}^{M_e \times M_e}, \quad [\mathbf{K}_e]_{ij} = \mathcal{L}[K(\mathbf{x}_i^{(e)}, \mathbf{x}_j^{(e)})], \\
\mathbf{S}_e &\in \mathbb{R}^{M_e \times |\mathcal{I}_e|}, \quad [\mathbf{S}_e]_{ij} = K(\mathbf{x}_i^{(e)}, \mathbf{x}_j^{\text{boundary}}), \\
\boldsymbol{B}_e &\in \mathbb{R}^{|\mathcal{I}_e| \times |\mathcal{I}_e|}, \quad [\boldsymbol{B}_e]_{ij} = K(\mathbf{x}_i^{\text{boundary}}, \mathbf{x}_j^{\text{boundary}}), \\
\mathbf{f}_e &\in \mathbb{R}^{M_e}, \quad [\mathbf{f}_e]_i = f(\mathbf{x}_i^{(e)}), \\
\mathbf{u}_e^{\text{fem}} &\in \mathbb{R}^{|\mathcal{I}_e|}, \quad [\mathbf{u}_e^{\text{fem}}]_i = u_i^{\text{fem}}, \\
\mathbf{b}_e &\in \mathbb{R}^{M_e}, \quad \mathbf{b}_e = [1, \ldots, 1]^T,
\end{align}
where $\mathcal{I}_e$ denotes the set of boundary node indices for element $e$, $|\mathcal{I}_e|$ is the number of boundary nodes per element, and $\mathbf{x}_j^{\text{boundary}}$ represents the boundary points of the element.
This system is solved for each element using parallel programming and KKT optimality conditions. This yields the Lagrangian multipliers $\boldsymbol{\alpha}_e \in \mathbb{R}^{M_e}$, $\boldsymbol{\beta}_e \in \mathbb{R}^{|\mathcal{I}_e|}$, and bias term $b_e \in \mathbb{R}$ for each element.
Upon solving the dual system for each element, the local enhanced solution within element $\Omega_e$ is given by
\begin{equation}
u_e(\mathbf{x}) = \sum_{j=1}^{M_e} \alpha_{e,j} K(\mathbf{x}, \mathbf{x}_j^{(e)}) + \sum_{k \in \mathcal{I}_e} \beta_{e,k} K(\mathbf{x}, \mathbf{x}_k) + b_e, \quad \mathbf{x} \in \Omega_e.
\label{eq:element_model}
\end{equation}

Step 3: After applying the method for each element, the hybrid global solution is defined as
\begin{equation}
u^{\text{hybrid}}(\mathbf{x}) = \begin{cases}
u_e(\mathbf{x}) & \text{if } \mathbf{x} \in \Omega_e \\
\frac{1}{|\mathcal{E}(\mathbf{x})|} \sum_{e \in \mathcal{E}(\mathbf{x})} u_e(\mathbf{x}) & \text{if } \mathbf{x} \in \partial\Omega_e \text{ for multiple elements}
\end{cases}
\label{eq:global_construction}
\end{equation}
where $\mathcal{E}(\mathbf{x})$ is the set of elements that share the point $\mathbf{x}$. 
%
%
For implementation simplicity, we use NumPy's tools for Legendre polynomials and their derivatives. The full algorithm for this hybrid method is shown in Algorithm~\ref{alg:fem-lssvr}.

\begin{algorithm*}[t]
\caption{Hybrid FEM-LSSVR Method}
\label{alg:fem-lssvr}
\begin{algorithmic}[1]
\State Initialize: domain $\Omega$, PDE parameters, regularization parameter $\gamma$
\State Apply standard FEM to obtain nodal solution: $\{(\mathbf{x}_i^{\mathrm{fem}}, u_i^{\mathrm{fem}})\}_{i=1}^{n_{\mathrm{fem}}}$
\State Partition domain into finite elements: $\Omega = \bigcup_{e=1}^{E_{\mathrm{fem}}} \Omega_e$
\For{$e = 1$ to $E_{\mathrm{fem}}$ in parallel}
    \State Extract element nodes: $\mathcal{N}_e = \{\mathbf{x}_k : \mathbf{x}_k \in \partial\Omega_e \cap \{\mathbf{x}_i^{\mathrm{fem}}\}\}$
    \State Generate interior collocation points: $\{\mathbf{x}_j^{(e)}\}_{j=1}^{M_e} \subset \Omega_e$
    \State Formulate LSSVR optimization problem with Legendre polynomial kernel:
    \State \quad $\min_{\mathbf{w},b,\mathbf{e}} \frac{1}{2}\mathbf{w}^T\mathbf{w} + \frac{\gamma}{2}\sum_{j=1}^{M} e_j^2$
    \State \quad s.t. PDE residual constraints at collocation points
    \State \quad and PDE boundary conditions on $\partial\Omega_\mathbf{e} \cap \partial\Omega$
    \State \quad and FEM nodal constraints at interior element boundaries
    \State Assemble element dual system from equation (\ref{eq:dual_system})
    \State Solve for Lagrange multipliers: $(\boldsymbol{\alpha}_e, \boldsymbol{\beta}_e, b_e)$
    \State Store element model: $u_e(\mathbf{x}) = \sum_j \alpha_{e,j} K(\mathbf{x}, \mathbf{x}_j^{(e)}) + \sum_k \beta_{e,k} K(\mathbf{x}, \mathbf{x}_k) + b_e$ from equation (\ref{eq:element_model})
\EndFor
\State Construct global solution function $u^{\mathrm{hybrid}}(\mathbf{x})$ using equation (\ref{eq:global_construction})
\State \Return Enhanced global solution $u^{\mathrm{hybrid}}(\mathbf{x})$
\end{algorithmic}
\end{algorithm*}
\subsection{Experimental Setup and Error Metrics}

Training and test points are independently selected locations where the solution accuracy is evaluated to measure the quality of training and generalization. For FEM, set called "training points" correspond to the finite element degrees of freedom (nodal locations). For LSSVR, training points are selected interior collocation points where the PDE constraints are enforced, using uniform distributions within the domain. For the hybrid FEM-LSSVR method, training points include both the FEM nodal values (which serve as boundary constraints for the element-wise LSSVR enhancement) and the interior collocation points within each finite element subdomain. Test points are uniformly distributed using linearly spaced grids across the computational domain, with the number of test points specified for each experiment. We also record execution time to measure the computational efficiency of each method.
The kernel order and regularization parameter $\gamma$ are selected through parameter optimization using Optuna \cite{optuna}, which employs a tree-structured Parzen estimator algorithm to efficiently explore the hyperparameter space and identify optimal configurations for each test case.

For each example, we evaluate performance using both relative $L_2$ and $H^1$ error metrics. The relative $L_2$ error 
\begin{equation}
\varepsilon_{\text{rel}}^{L^2} = \frac{\|u_h - u\|_{L^2(\Omega)}}{\|u\|_{L^2(\Omega)}} = \frac{\sqrt{\int_{\Omega} \left(u_h(x) - u(x)\right)^2 \, d\Omega}}{\sqrt{\int_{\Omega} u(x)^2 \, d\Omega}},
\end{equation}
serves as the standard convergence metric in finite element analysis and serves as an objective function in LSSVR optimization.
We evaluate the $L_2$ norm by numerical integration using \texttt{scipy.integrate.quad} where $u_h(x)$ is the numerical solution (from FEM, LSSVR, or the hybrid FEM-LSSVR method), $u(x)$ is the exact analytical solution. 
The relative $H^1$ error 
\begin{equation}
\varepsilon_{\text{rel}}^{H^1} = \frac{\|u_h - u\|_{H^1(\Omega)}}{\|u\|_{H^1(\Omega)}} = \frac{\sqrt{\int_{\Omega} \left[(u_h(x) - u(x))^2 + |\nabla u_h(x) - \nabla u(x)|^2\right] \, d\Omega}}{\sqrt{\int_{\Omega} \left[u(x)^2 + |\nabla u(x)|^2\right] \, d\Omega}},
\end{equation}
provides additional insight into the accuracy of the solution gradients. The metrics outcome should demonstrate that the hybrid method maintains or improves upon the theoretical convergence rates of both FEM and LSSVR across all tested PDE examples.


\section{Numerical Results}
\label{sec:results}

To examine the performance of the hybrid FEM-LSSVR method, we conducted numerical experiments on four PDE examples, including Poisson equations in 1D and 2D domains, as well as Helmholtz equations with various boundary conditions. The implementation consumed 64 cores, 4 GB of RAM, and a benchmark value (Rmax) of 2.31 PFlop/s. In this section, we present results of the hybrid approach with a standard (low-order) FEM implementation and the standard LSSVR method.

\subsection{1D Poisson equation}

We first consider the 1D Poisson equation, 
\begin{equation*}
    -\frac{\mathrm{d}^2u}{\mathrm{d}x^2} = f(x), \quad x \in [a,b]
    \, ,
\end{equation*}
with $f(x) = \pi^2\sin{(\pi x)}$ and 
%
%
with homogeneous Dirichlet boundary conditions
$
    u(a) = 0, \, u(b) = 0
    .
$
The exact solution to this problem is
$
    u(x) = \sin(\pi x)
    .
$
We solved this problem in two different domains, $[-1,1]$ and $[-5,5]$, to test performance over a large domain. For the FEM implementation, we used linear (P1), quadratic (P2), and fifth-order (P5) elements to compare methods with the same polynomial order. For the LSSVR method, we used 5 collocation points, and for the hybrid FEM+LSSVR method, we used 5 collocation with $\gamma = 10^{6}$ points in each element. Table~\ref{tab1} shows performance results for domains $[-1,1]$ and $[-5,5]$. In both cases, we used 24 elements for the FEM and the hybrid FEM+LSSVR methods.

\begin{figure}[ht!]
\centering
\begin{subfigure}[b]{0.3\textwidth}
    \includegraphics[width=\textwidth]{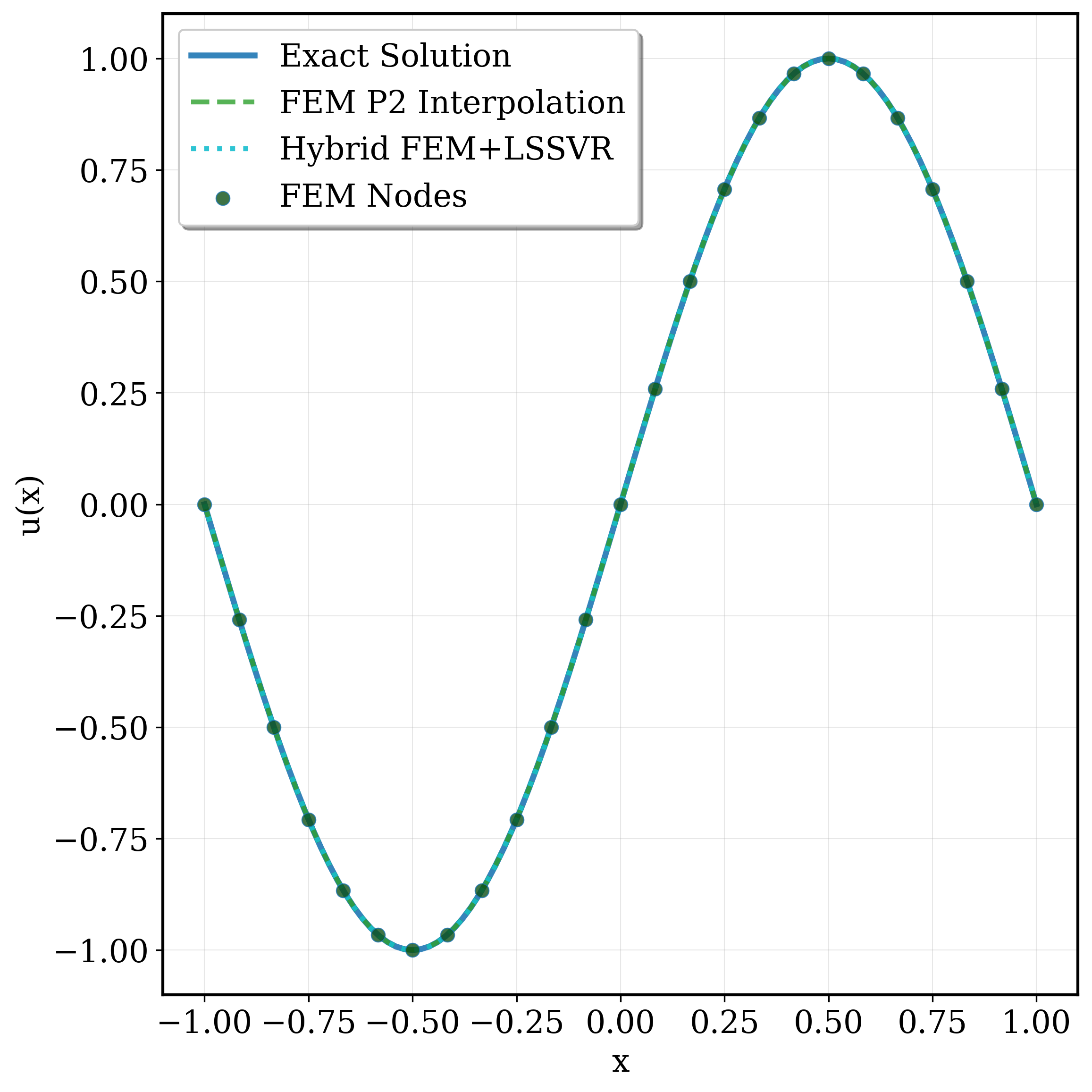}
    \caption{Fine Mesh Solution on [-1,1]}
    \label{fig:ex1_hybrid_sol_1}
\end{subfigure}
\hfill
\begin{subfigure}[b]{0.3\textwidth}
    \includegraphics[width=\textwidth]{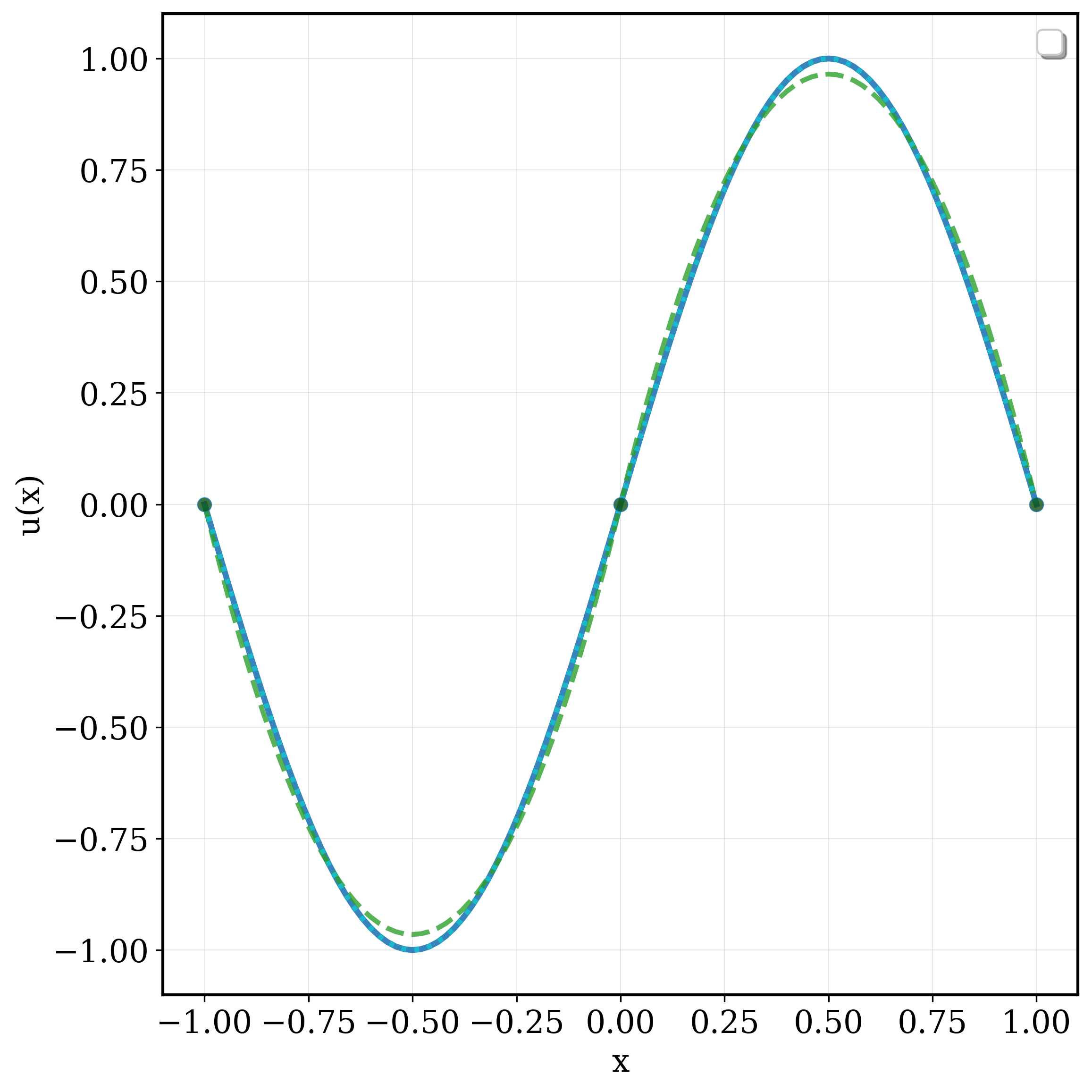}
    \caption{Coarse Mesh Solution on [-1,1]}
    \label{fig:ex1_hybrid_err_1}
\end{subfigure}
\hfill
\begin{subfigure}[b]{0.3\textwidth}
    \includegraphics[width=\textwidth]{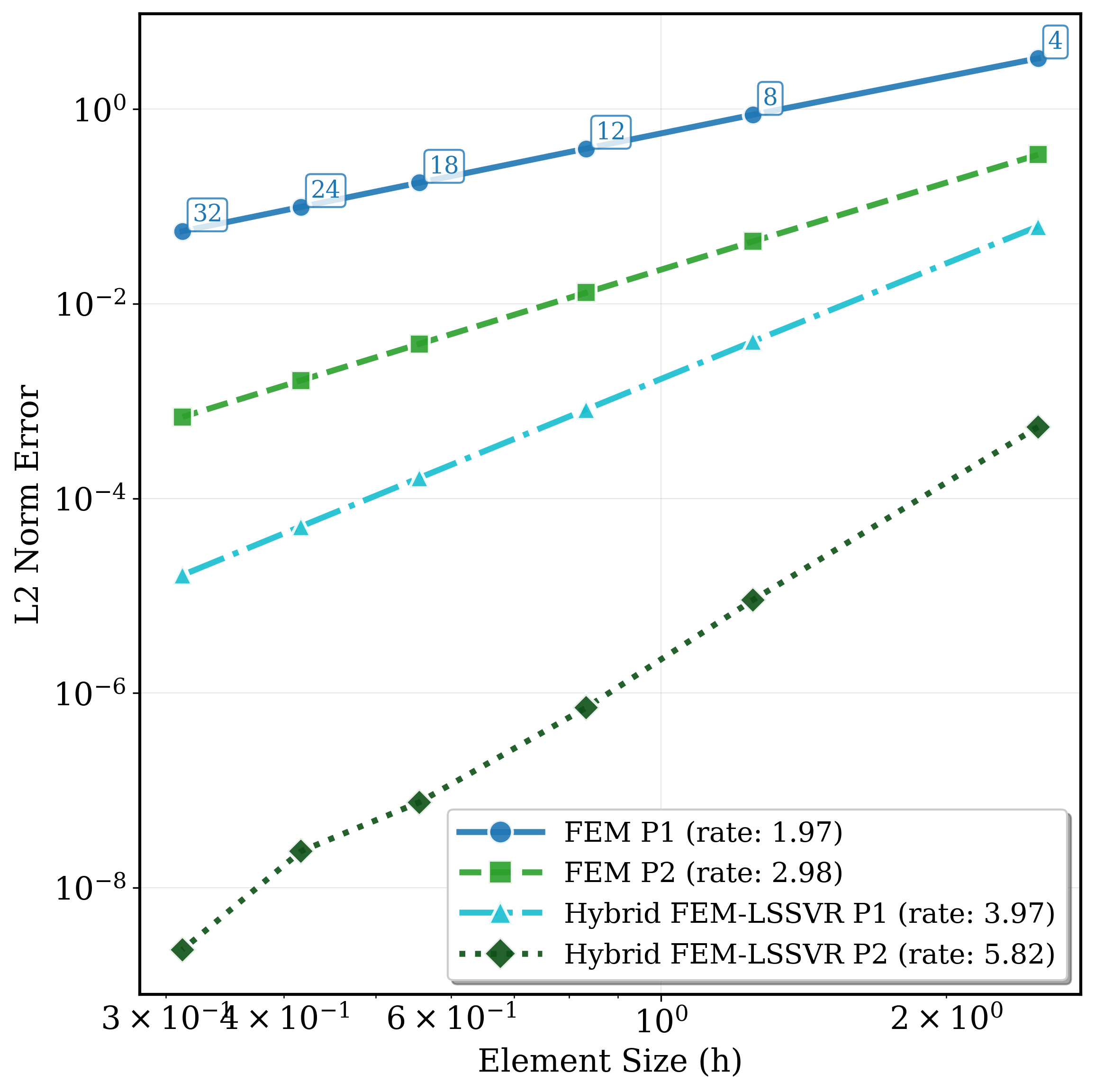}
    \caption{Convergence Analysis on [-1,1]}
    \label{fig:ex1_convergence_1}
\end{subfigure}

\vspace{0.5cm} 

\begin{subfigure}[b]{0.3\textwidth}
    \includegraphics[width=\textwidth]{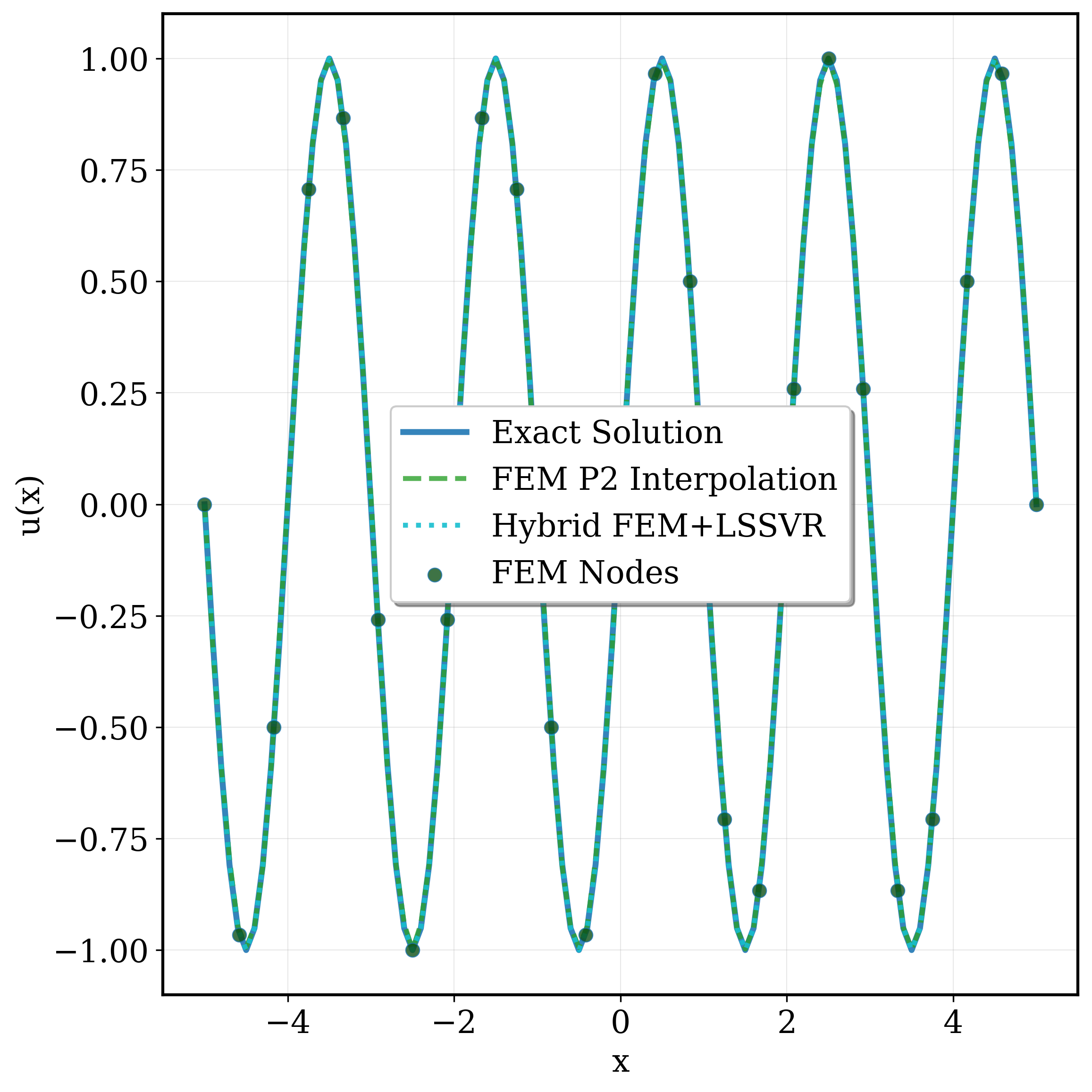}
    \caption{Fine Mesh Solution on [-5,5]}
    \label{fig:ex1_hybrid_sol_5}
\end{subfigure}
\hfill
\begin{subfigure}[b]{0.3\textwidth}
    \includegraphics[width=\textwidth]{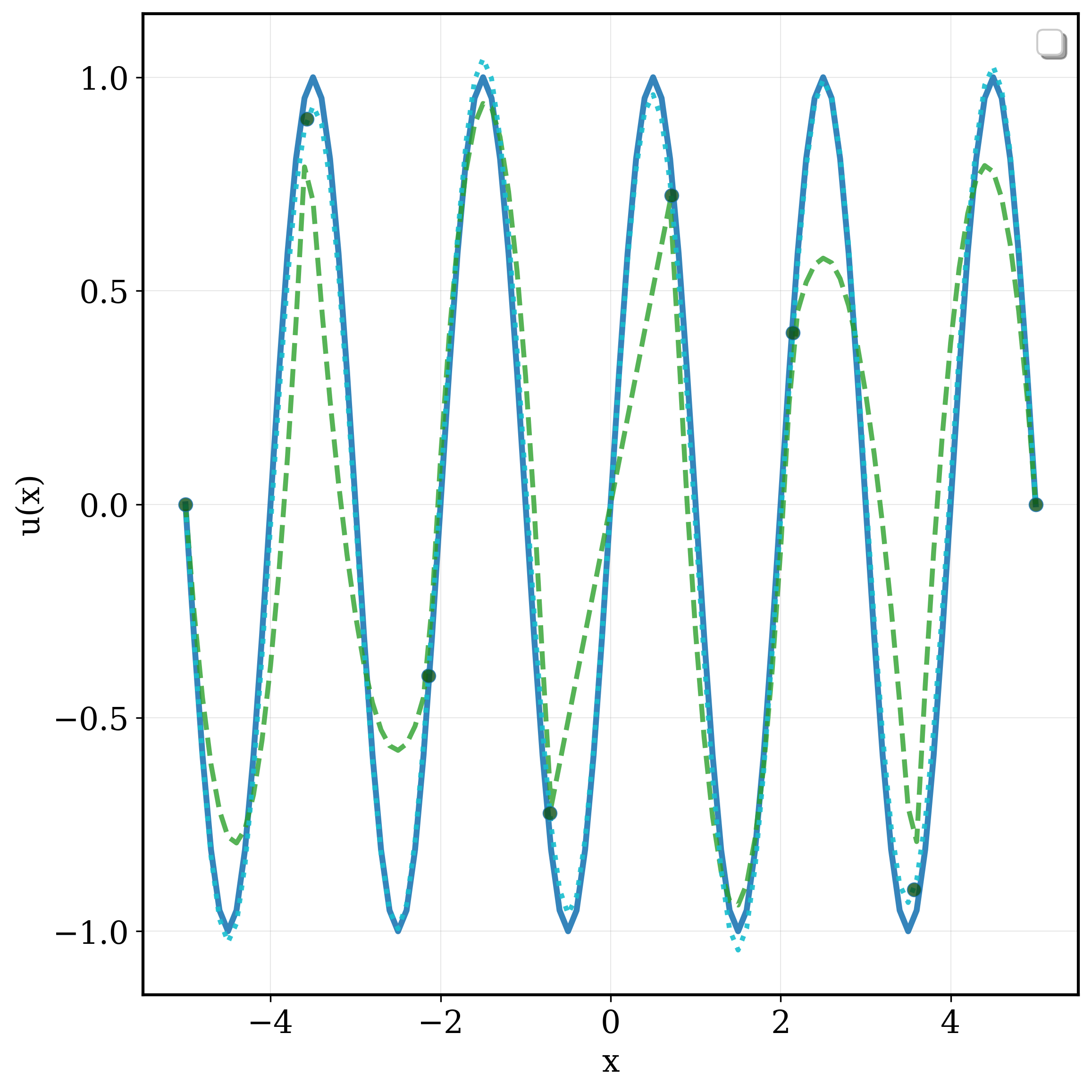}
    \caption{Coarse Mesh Solution on [-5,5]}
    \label{fig:ex1_hybrid_err_5}
\end{subfigure}
\hfill
\begin{subfigure}[b]{0.3\textwidth}
    \includegraphics[width=\textwidth]{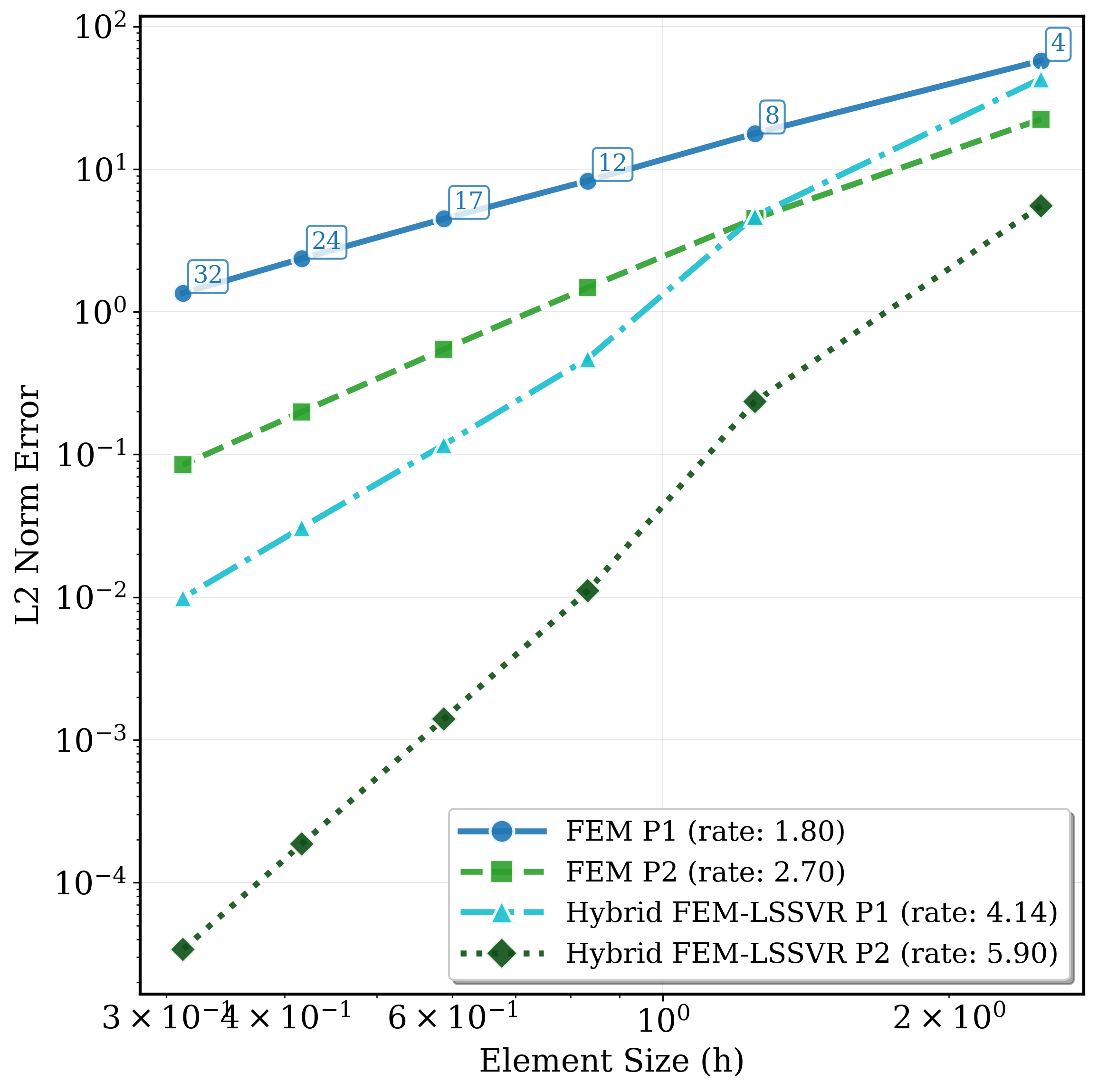}
    \caption{Convergence Analysis on [-5,5]}
    \label{fig:ex1_convergence_5}
\end{subfigure}

\caption{Example 1 - 1D Poisson Equation Results: Comparison of fine mesh (many FEM nodes), coarse mesh (fewer FEM nodes), and convergence analysis for domains [-1,1] (top row) and [-5,5] (bottom row)}

\label{fig:ex1_results}
\end{figure}

Figure~\ref{fig:ex1_results} visualizes the computational accuracy of the hybrid FEM-LSSVR approach across various mesh densities and domains. On finer meshes, both the FEM with second-order (P2) interpolation and the hybrid method perform well and captured the underlying sinusoidal pattern. The hybrid approach demonstrates improved accuracy compared to standalone FEM and LSSVR implementations. Additionally, the results on coarse meshes proved the super-resolution advantage of our method. A standard (low-order) FEM normally requires dense meshes to achieve accurate results, whereas the hybrid method provides accurate solutions even with a minimal number of elements. This is particularly evident in the $[-5,5]$ domain, where the coarse FEM solution with a few nodes exhibits notable oscillation errors and a poor representation of the underlying function. In contrast, the hybrid solution with the same sparse mesh maintains visibly accurate results.
In contrast, the hybrid method keeps almost perfect accuracy even with sparse data.
Based on the convergence analysis  while FEM with P1 and P2 elements shows the expected polynomial convergence rates ($O(h^2)$ and $O(h^3)$, respectively), the hybrid method achieves improved convergence across both domains, maintaining accuracy that is 3–6 orders of magnitude higher than standard (low-order) FEM across all levels of mesh refinement.


Table~\ref{tab1} shows results for both domains. In the smaller domain, the hybrid methods achieve 3--4 orders of magnitude better accuracy than standard low-order FEM (P1 and P2) across all relative $L_2$ error measures, with FEM P2+LSSVR performing the best. In addition, FEM P2+LSSVR outperforms FEM P5 via the fifth-order Legendre polynomial kernels, showing the superior interpolation capabilities of these kernels. 
When the domain is expanded to $[-5,5]$, the relative errors generally increase. Standard FEM P1 and P2 lose significant accuracy, with errors rising by several orders of magnitude, while the hybrid methods continue to provide much more accurate results. The standalone LSSVR performed very poorly for the larger domain. However, the hybrid methods maintain high accuracy in both the $L^2$ and $H^1$ norms, demonstrating smooth and accurate solution gradients. The computation times to obtain the solution were comparable over all methods for the larger domain. We experienced faster evaluation times for the LSSVR and hybrid models compared to FEM (which uses scikit-fem with a direct solver). On the small domain, the FEM simulation was significantly faster than the LSSVR and hybrid methods in obtaining the solution.
 
\begin{table*}[ht!]
\centering
\caption{Performance comparison of FEM, LSSVR, and hybrid FEM+LSSVR methods for solving the 1D Poisson equation.\label{tab1}}
\begin{tabular*}{\textwidth}{@{\extracolsep\fill}lcccccc@{\extracolsep\fill}}
\toprule
\textbf{Description/Metric} & \textbf{FEM P1} & \textbf{FEM P2} & \textbf{FEM P5} & \textbf{LSSVR} & \textbf{FEM P1+LSSVR} & \textbf{FEM P2+LSSVR} \\
\midrule
Domain [-1,1] \\
\midrule
Basis Order & 1 & 2 & 5 & -- & 1 & 2 \\
Kernel Order & -- & -- & -- & 5 & 5 & 5 \\
\midrule
Computation Time (s) & 0.001771 & 0.001899 & 0.003321 & 0.103894 & 0.286176 & 0.295793 \\
Evaluation Time (s)  & 0.021643 & 0.026769 & 0.045824 & 0.000068 & 0.006755 & 0.006716 \\
\midrule
Solution $\varepsilon_{\text{rel}}$ & 6.244399$\times$10$^{-3}$ & 1.044729$\times$10$^{-4}$ & 1.964984$\times$10$^{-8}$ & 1.269698$\times$10$^{-3}$ & 3.255423$\times$10$^{-6}$ & 1.008775$\times$10$^{-9}$ \\
Relative H1 Error & 7.200588$\times$10$^{-2}$ & 2.436121$\times$10$^{-3}$ & 9.208956$\times$10$^{-7}$ & 2.803557$\times$10$^{-3}$ & 3.263940$\times$10$^{-6}$ & 2.080386$\times$10$^{-8}$ \\
\midrule
Domain [-5,5] \\
\midrule
Computation Time (s) & 0.001729 & 0.001837 & 0.003210 & 0.077334 & 0.189902 & 0.216317 \\
Evaluation Time (s)  & 0.023751 & 0.025684 & 0.046998 & 0.981808 & 0.007053 & 0.006830 \\
\midrule
Solution $\varepsilon_{\text{rel}}$ & 1.489940$\times$10$^{-1}$ & 1.257327$\times$10$^{-2}$ & 6.042536$\times$10$^{-5}$ & 2.184097$\times$10$^{2}$ & 1.940345$\times$10$^{-3}$ & 1.525231$\times$10$^{-5}$ \\
Relative H1 Error & 3.525858$\times$10$^{-1}$ & 5.907133$\times$10$^{-2}$ & 7.876735$\times$10$^{-4}$ & 8.486228$\times$10$^{1}$ & 2.068411$\times$10$^{-3}$ & 6.524946$\times$10$^{-5}$ \\
\bottomrule
\end{tabular*}
\end{table*}

\subsection{2D Poisson equation}
As a second example, we next consider the 2D Poisson equation \eqref{eq:poisson}. 
%
%
Specifically, we solve
\begin{equation*}
    -\nabla^2 u=2\pi^2\sin{(\pi x)}\sin{(\pi y)}, \quad (x,y) \in [a,b] \times [c,d] ,
\end{equation*}
with Dirichlet boundary conditions
%
    $u(a,y)=u(b,y)=u(x,c)=u(x,d)=0$, $ \forall x,y$. 
%
The exact solution is
%
    $u(x,y)=\sin{(\pi x)}\sin{(\pi y)}$.
%
This equation is solved in two domains $[-1,1] \times [-1,1]$ and $[-2,2] \times [-2,2]$ to test how the methods handle increased domain size. For the LSSVR method, we used $6 \times 6$ uniformly distributed collocation points with $\gamma = 10^{5}$, and for the hybrid method, we used the same number in each element. 

\begin{figure}[ht!]
\centering
\begin{subfigure}[b]{0.32\textwidth}
    \includegraphics[width=\textwidth]{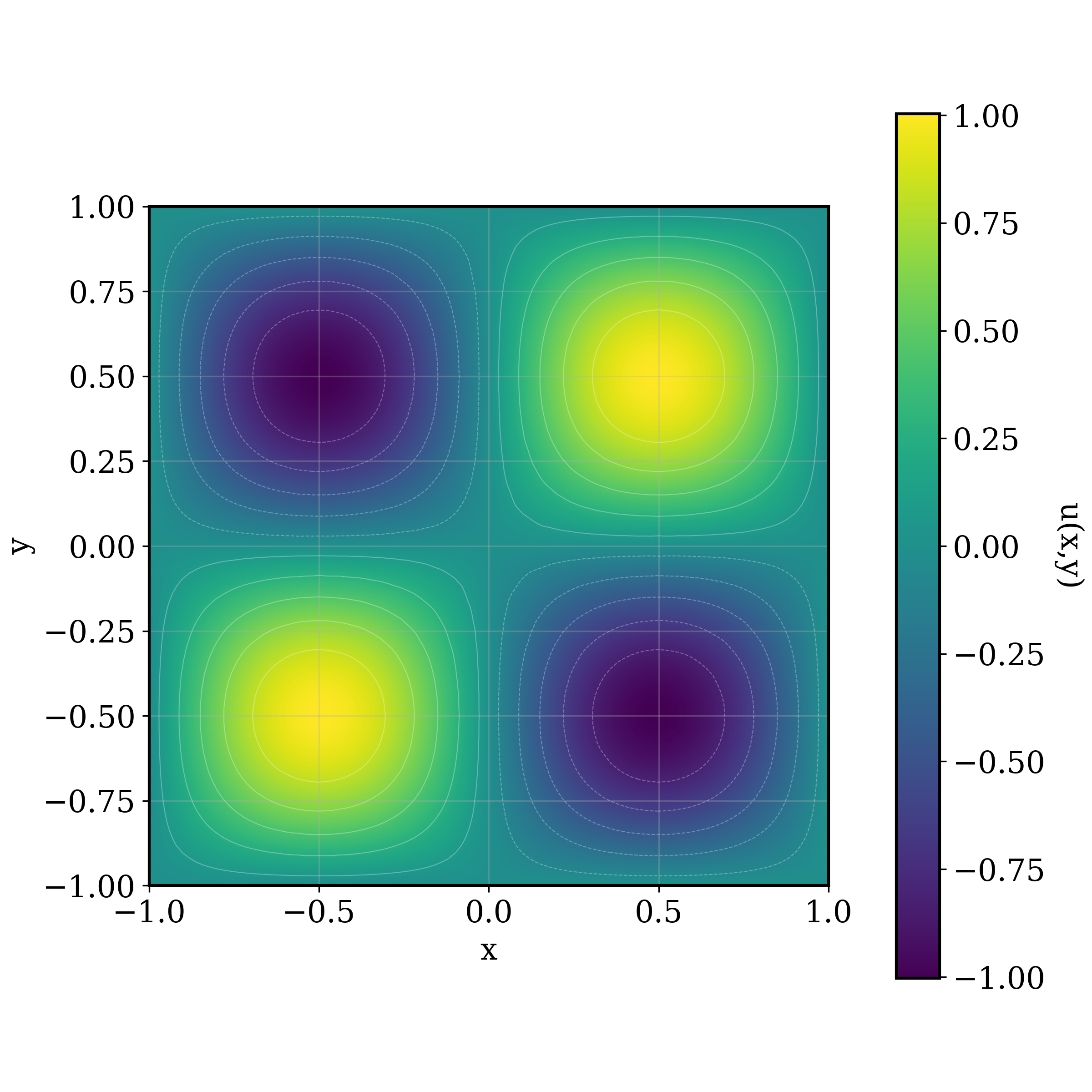}
    \caption{Hybrid Solution on [-1,1]$\times$[-1,1]}
    \label{fig:ex2_hybrid_sol_1}
\end{subfigure}
\hfill
\begin{subfigure}[b]{0.32\textwidth}
    \includegraphics[width=\textwidth]{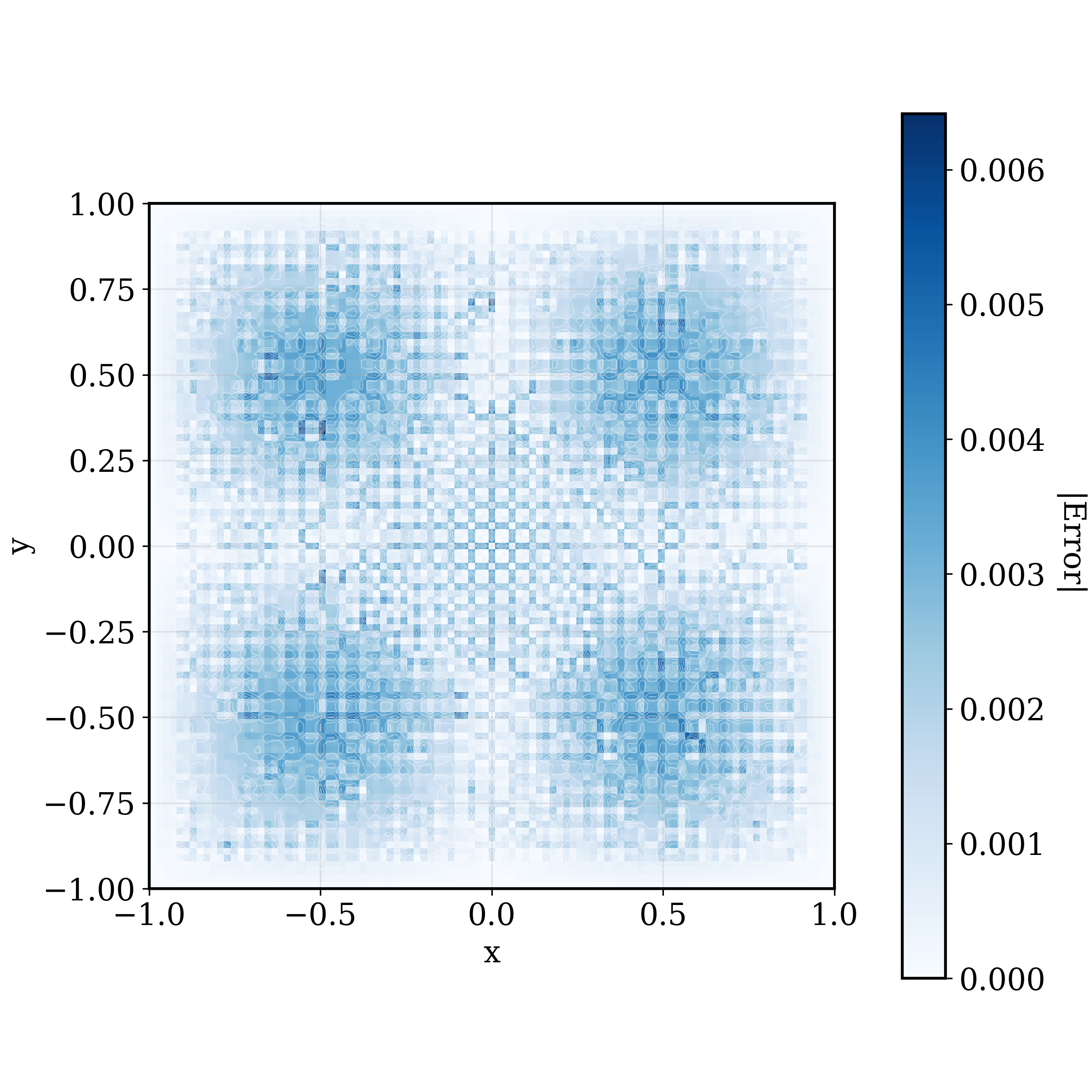}
    \caption{Error Field on [-1,1]$\times$[-1,1]}
    \label{fig:ex2_hybrid_err_1}
\end{subfigure}
\hfill
\begin{subfigure}[b]{0.32\textwidth}
    \includegraphics[width=\textwidth]{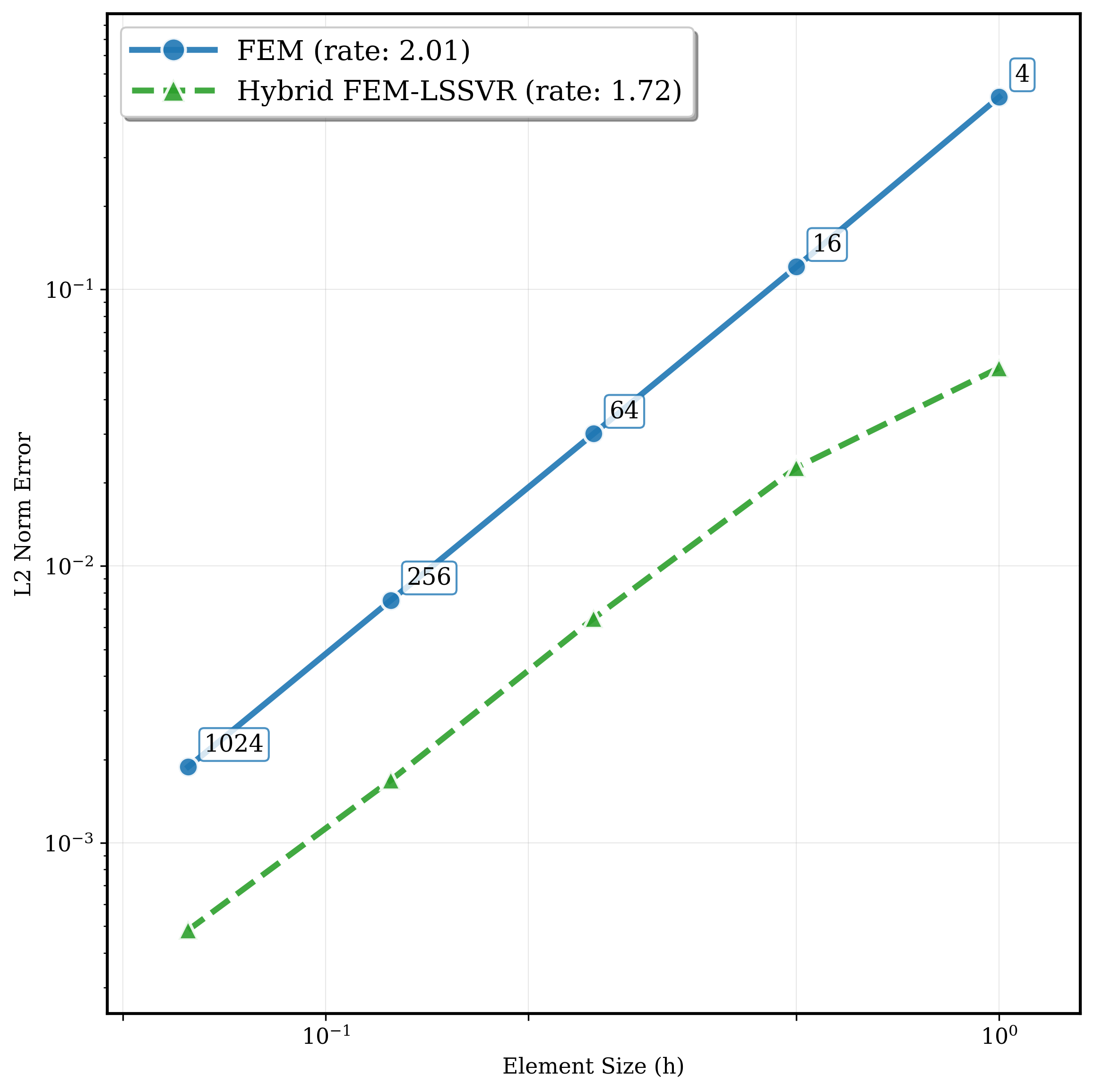}
    \caption{Convergence Analysis on [-1,1]$\times$[-1,1]}
    \label{fig:ex2_convergence_1}
\end{subfigure}
\vspace{0.5cm} 
\begin{subfigure}[b]{0.32\textwidth}
    \includegraphics[width=\textwidth]{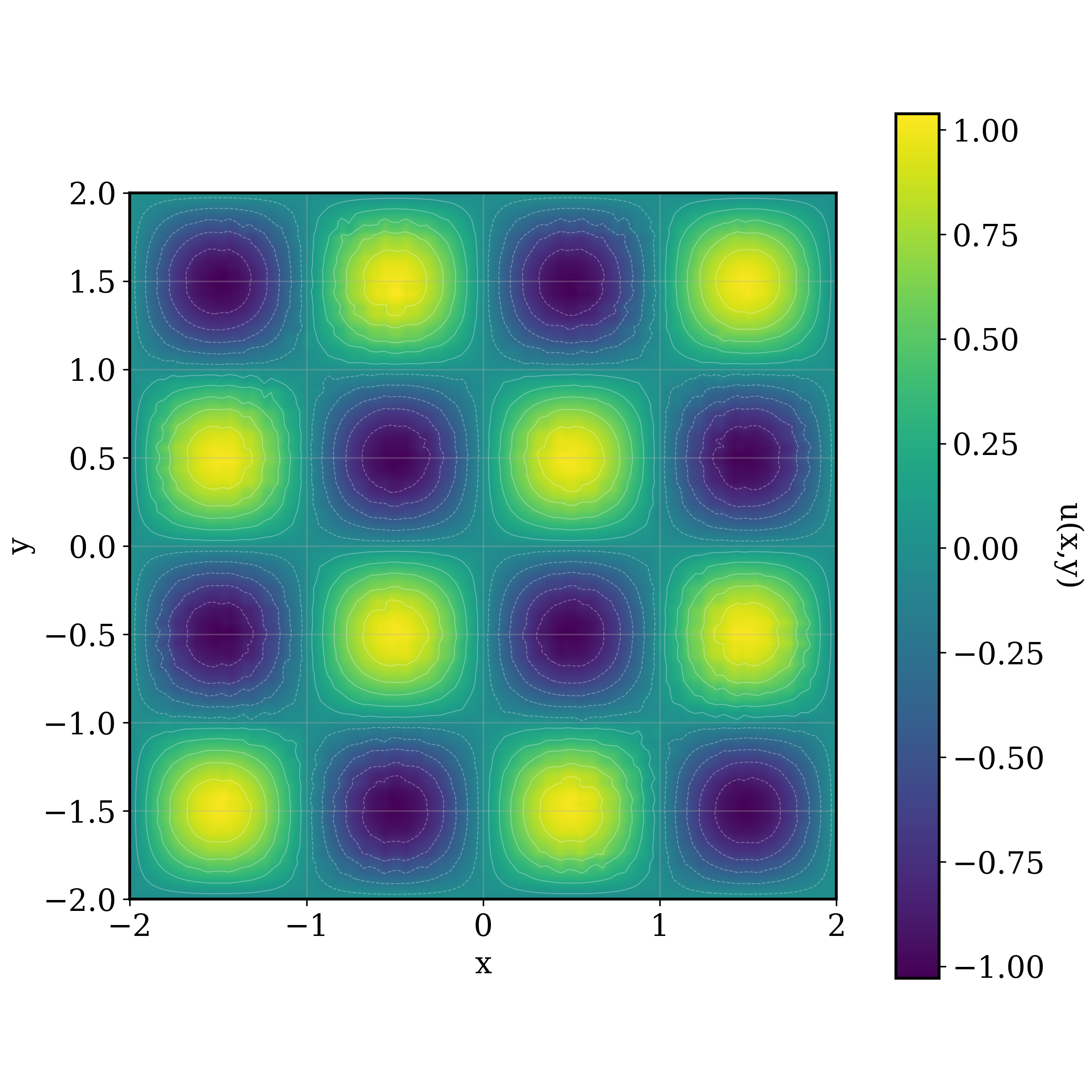}
    \caption{Hybrid Solution on [-2,2]$\times$[-2,2]}
    \label{fig:ex2_hybrid_sol_2}
\end{subfigure}
\hfill
\begin{subfigure}[b]{0.32\textwidth}
    \includegraphics[width=\textwidth]{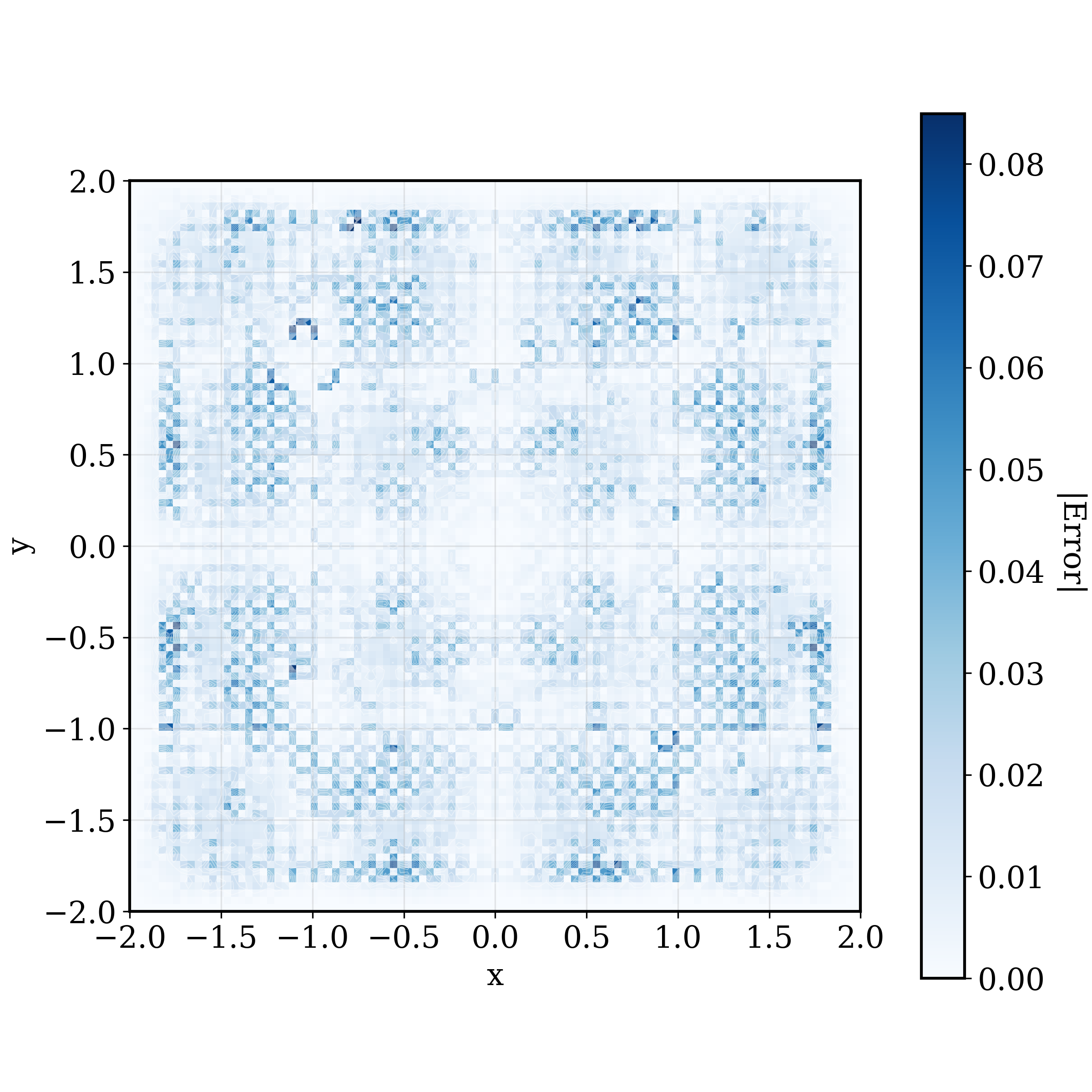}
    \caption{Error Field on [-2,2]$\times$[-2,2]}
    \label{fig:ex2_hybrid_err_2}
\end{subfigure}
\hfill
\begin{subfigure}[b]{0.32\textwidth}
    \includegraphics[width=\textwidth]{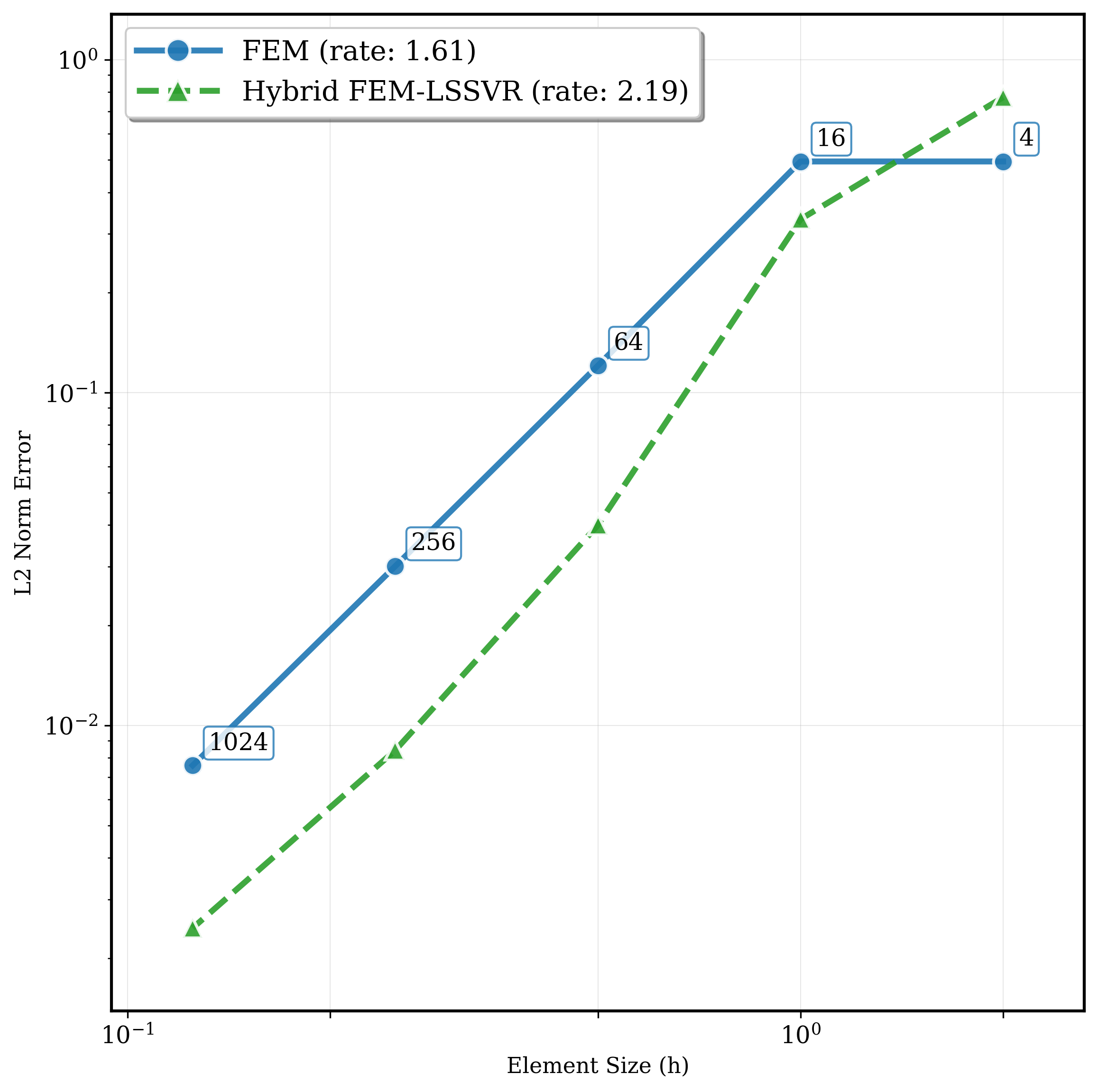}
    \caption{Convergence Analysis on [-2,2]$\times$[-2,2]}
    \label{fig:ex2_convergence_2}
\end{subfigure}

\caption{Example 2 – 2D Poisson Equation Results: Hybrid FEM-LSSVR solution, error field, and convergence analysis for domains [-1,1]$\times$[-1,1] (top row) and [-2,2]$\times$[-2,2] (bottom row).}
\label{fig:ex2_results}
\end{figure}

Figure~\ref{fig:ex2_results} visualizes the 2D solution, error distribution and convergence analysis for both domains. The hybrid method captures the sinusoidal pattern very well, with only tiny errors scattered lightly across the domain. The convergence analysis demonstrates the accuracy advantage over traditional FEM, showing steady improvement as resolution increases without any convergence stagnation.

Table~\ref{tab3} presents the results for the 2D domains using 1024 elements. For the smaller domain, the hybrid method achieves the best solution accuracy at $7.49 \times 10^{-4}$, outperforming FEM ($3.80 \times 10^{-3}$) by approximately 5 times and LSSVR ($3.40 \times 10^{-2}$) by about 45 times. For checking the smoothness of the method, we also consider the $H^1$ norm. The hybrid method achieves errors that are approximately 4 times better than FEM and 20 times better than standalone LSSVR. In addition to providing an accurate solution, the hybrid approach also provides high-accuracy spatial derivatives.

When the domain size increases, both the FEM and the separate LSSVR method perform poorly, with errors sharply increasing by orders of magnitude. The hybrid method achieves approximately two times the accuracy of FEM in the solution and retains better accuracy in the $H^1$ norm. In the larger domain, the individual LSSVR method shows limited accuracy, demonstrating that machine learning methods can be sensitive to domain size and model complexity.

\begin{table*}[ht!]
\centering
\caption{Performance comparison of FEM, LSSVR, and hybrid FEM+LSSVR methods for solving the 2D Poisson equation.\label{tab3}}
\begin{tabular*}{\textwidth}{@{\extracolsep\fill}lccc@{\extracolsep\fill}}
\toprule
\textbf{Description/Metric} & \textbf{FEM} & \textbf{LSSVR} & \textbf{FEM+LSSVR} \\
\midrule
Domain [-1,1]$\times$[-1,1] \\
\midrule
Basis Order & 1 & -- & 1 \\
Kernel Order & -- & 6 & 3 \\
\midrule
Computation Time (s) & 0.025962 & 0.003777 & 7.384207 \\
Evaluation Time (s)  & 0.266772 & 0.164447 & 2.404072 \\
\midrule
Solution $\varepsilon_{\text{rel}}$ & 3.801147$\times$10$^{-3}$ & 3.399075$\times$10$^{-2}$ & 7.492093$\times$10$^{-4}$ \\
Relative H1 Error & 1.053354$\times$10$^{-2}$ & 5.715443$\times$10$^{-2}$ & 2.788450$\times$10$^{-3}$ \\
\midrule
Domain [-2,2]$\times$[-2,2] \\
\midrule
Computation Time (s) & 0.045926 & 0.037071 & 9.357472 \\
Evaluation Time (s)  & 0.281154 & 0.087016 & 2.341327 \\
\midrule
Solution $\varepsilon_{\text{rel}}$ & 1.532992$\times$10$^{-2}$ & 8.226224$\times$10$^{-1}$ & 7.390370$\times$10$^{-3}$ \\
Relative H1 Error & 4.898270$\times$10$^{-2}$ & 8.257240$\times$10$^{-1}$ & 3.852657$\times$10$^{-2}$ \\
\bottomrule
\end{tabular*}
\end{table*}

\subsection{2D Helmholtz equation with Dirichlet boundary conditions}

Furthermore, we tested the methods on the 2D Helmholtz equation \eqref{eq:helmholtz}. The Helmholtz equation with parameter $n$ is given by
\begin{equation*}
    -\nabla^2u-\xi^2u=f(x,y), \quad \xi = 2\pi n
    .
\end{equation*}
The source function is defined as
\begin{equation*}
    f(x,y)=\xi^2\sin{(\xi x)}\sin{(\xi y)}, \quad(x,y)\in[0,1]\times[0,1]
    ,
\end{equation*}
and the homogeneous Dirichlet boundary conditions read
\begin{equation*}
    u(0,y)=u(1,y)=u(x,0)=u(x,1)=0, \quad\forall x,y 
    .
\end{equation*}
The exact solution is
\begin{equation*}
    u(x,y)=\sin{(\xi x)}\sin{(\xi y)} .
\end{equation*}
This equation is solved in the domain $[0,1] \times [0,1]$ for two different wave numbers of $\xi$. These values correspond to different wavelengths in the solution, with higher values of n resulting in more oscillations. For the LSSVR method, we used $6 \times 6$ uniformly distributed collocation points with $\gamma = 4\times 10^{-4}$, and for the hybrid method, we used the same number in each element..

\begin{figure}[ht!]
\centering
\begin{subfigure}[b]{0.32\textwidth}
    \includegraphics[width=\textwidth]{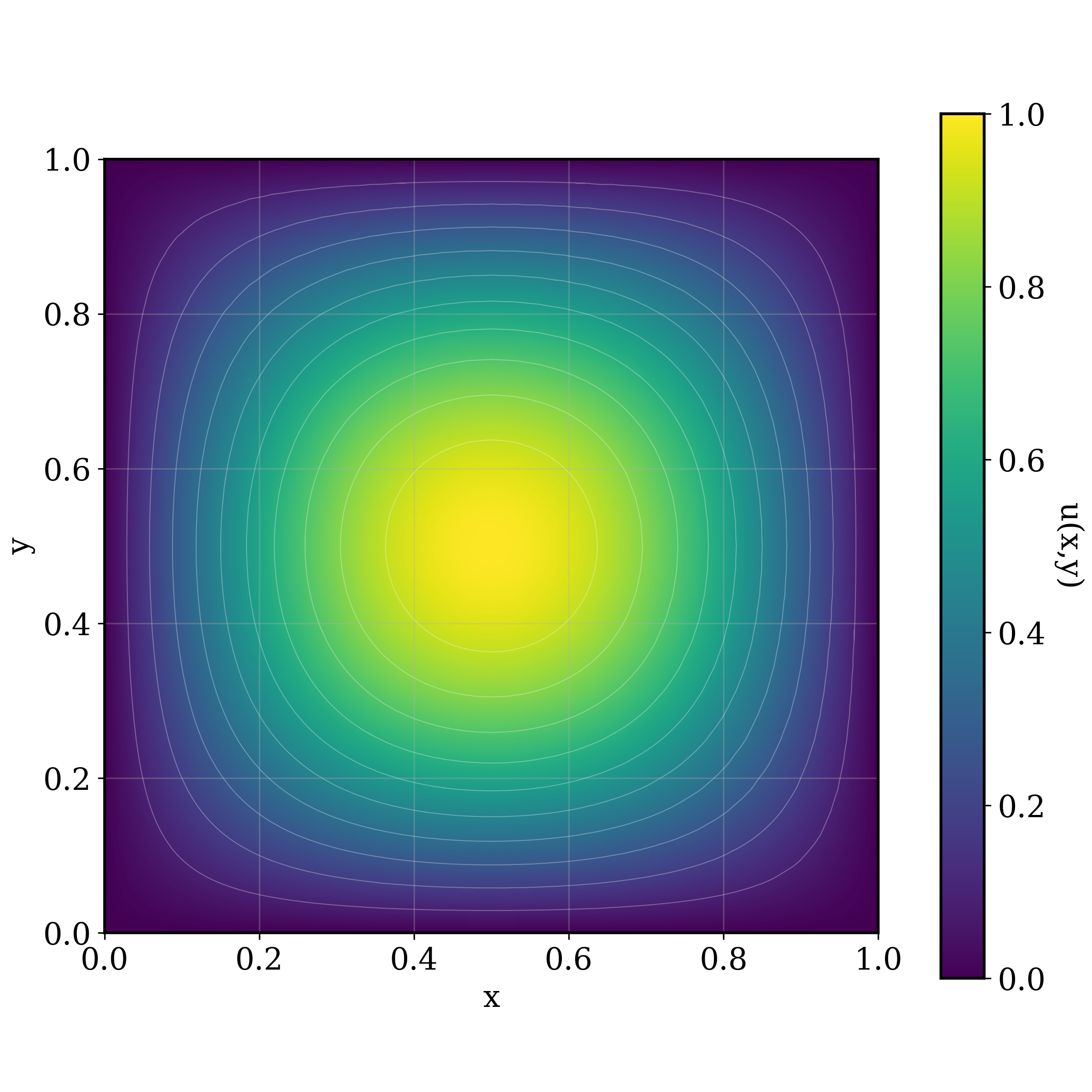}
    \caption{Hybrid Solution ($n = \frac{1}{2}$)}
    \label{fig:ex3_hybrid_sol_half}
\end{subfigure}
\hfill
\begin{subfigure}[b]{0.32\textwidth}
    \includegraphics[width=\textwidth]{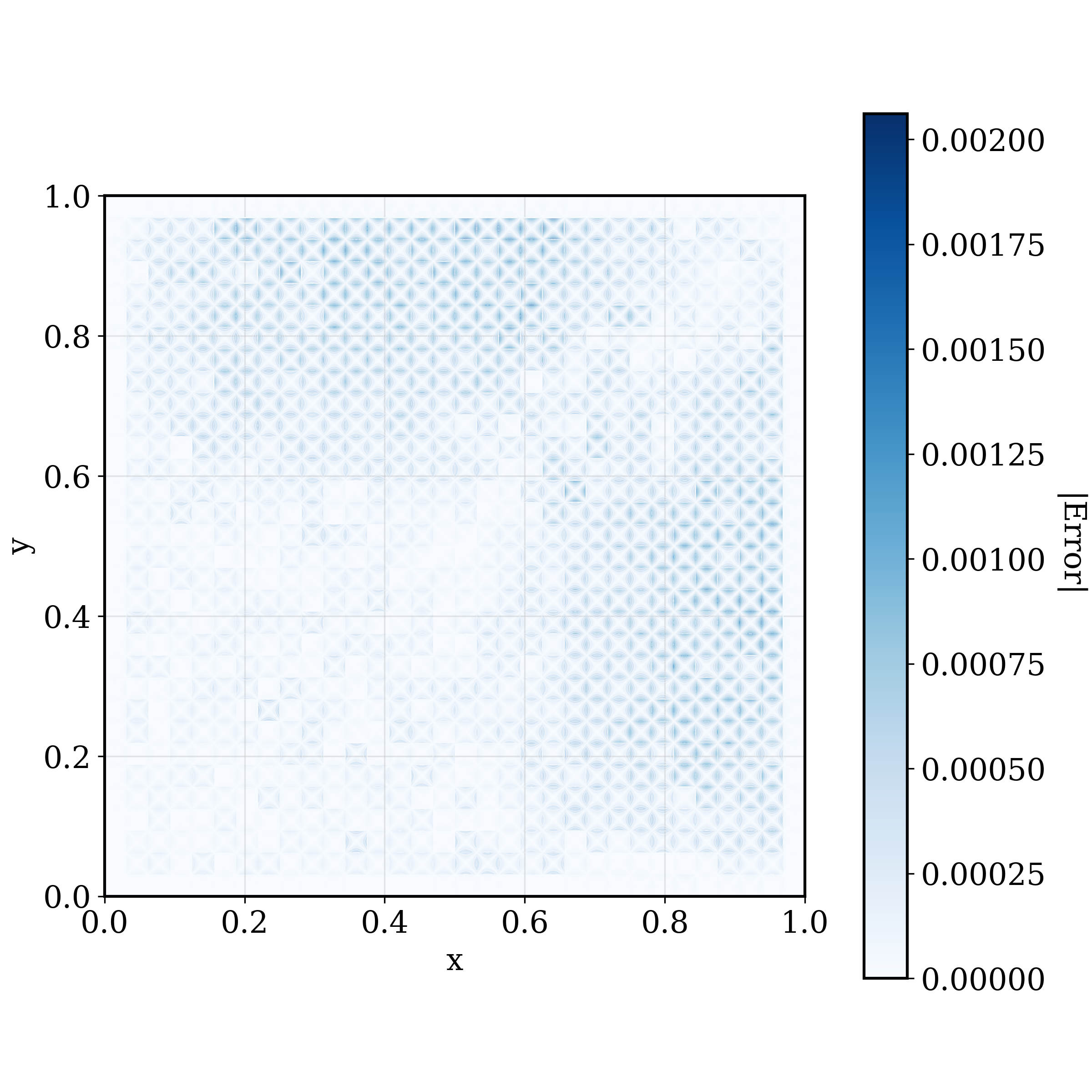}
    \caption{Error Field ($n = \frac{1}{2}$)}
    \label{fig:ex3_hybrid_err_half}
\end{subfigure}
\hfill
\begin{subfigure}[b]{0.32\textwidth}
    \includegraphics[width=\textwidth]{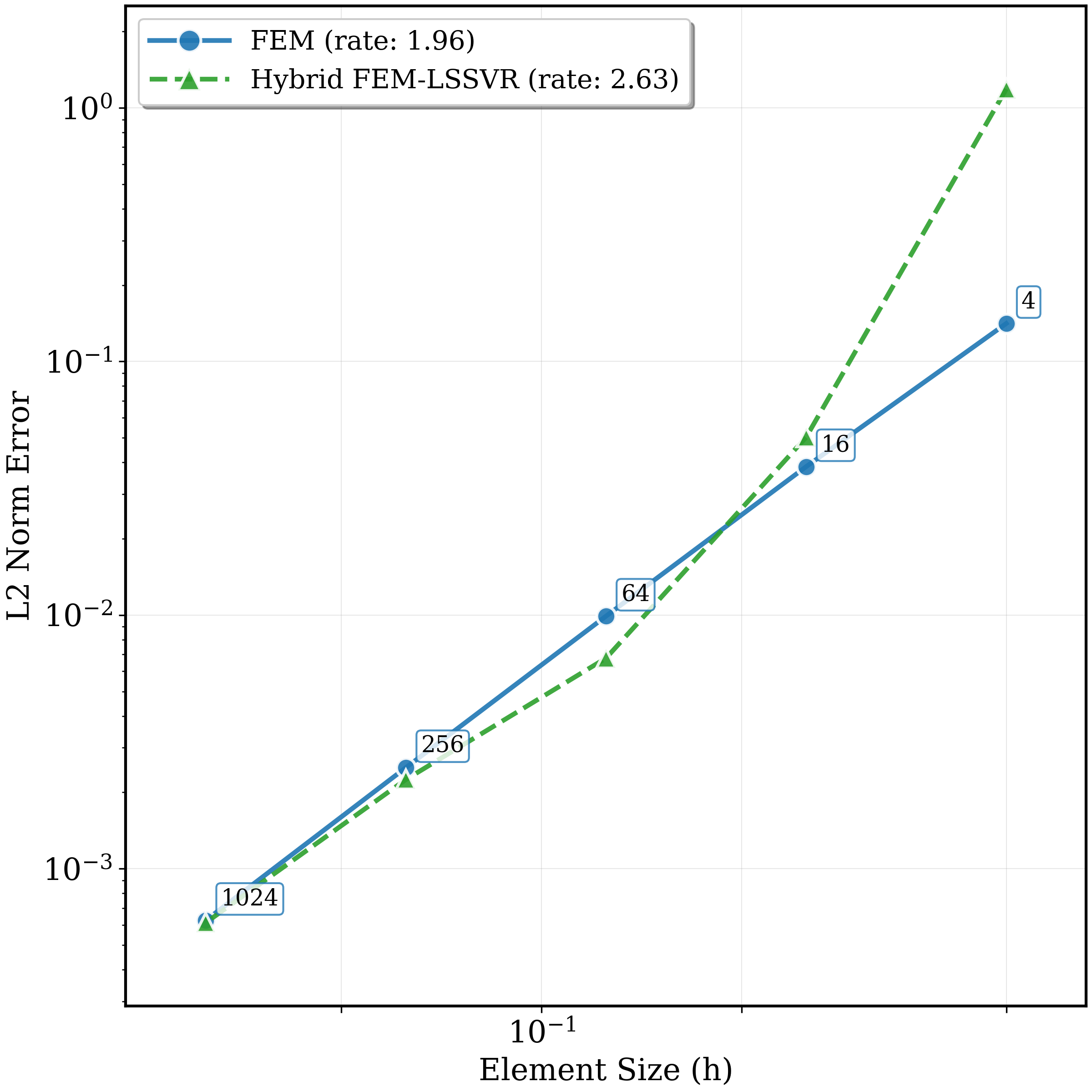}
    \caption{Convergence Analysis ($n = \frac{1}{2}$)}
    \label{fig:ex3_convergence_half}
\end{subfigure}

\vspace{0.5cm} 

\begin{subfigure}[b]{0.32\textwidth}
    \includegraphics[width=\textwidth]{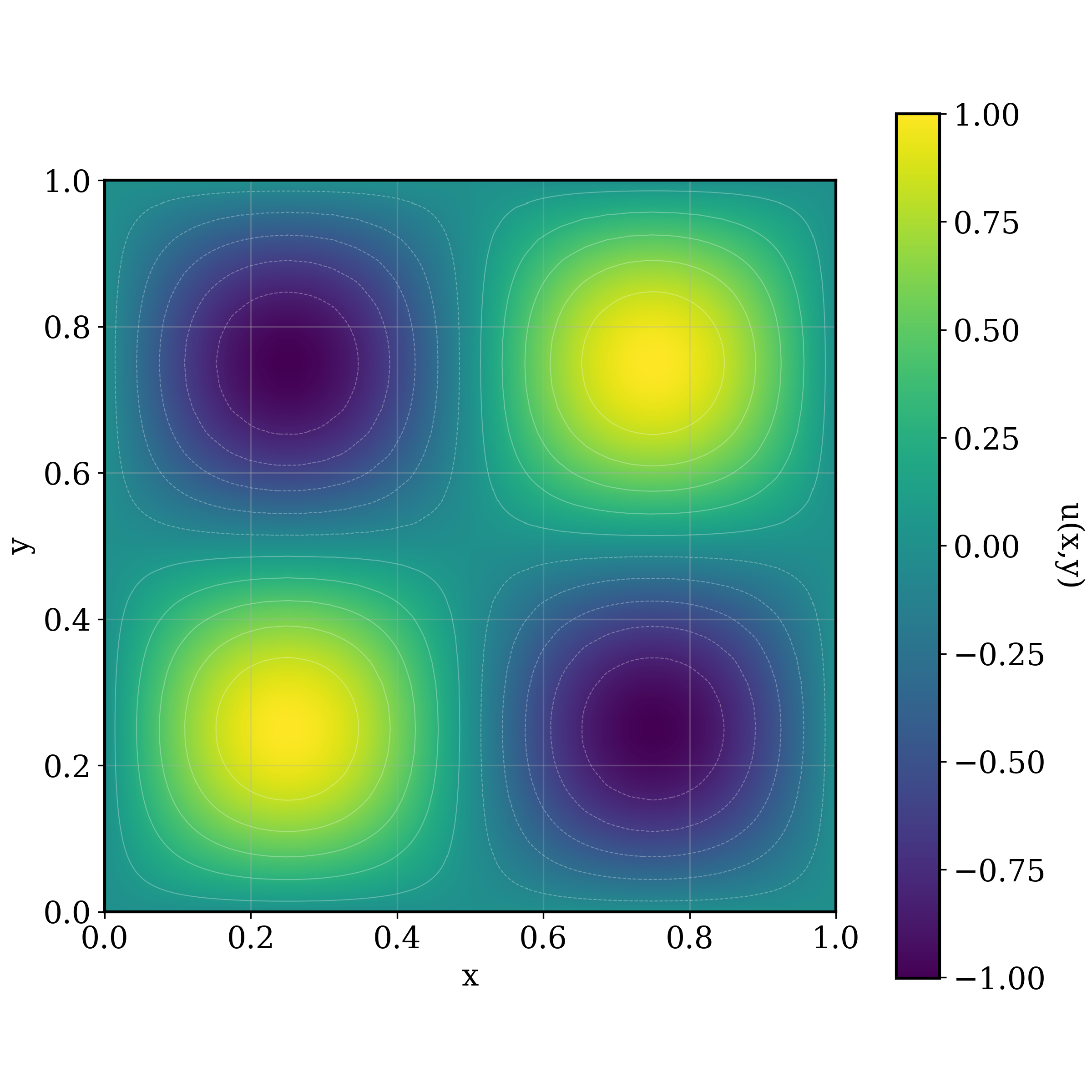}
    \caption{Hybrid Solution ($n = 1$)}
    \label{fig:ex3_hybrid_sol_one}
\end{subfigure}
\hfill
\begin{subfigure}[b]{0.32\textwidth}
    \includegraphics[width=\textwidth]{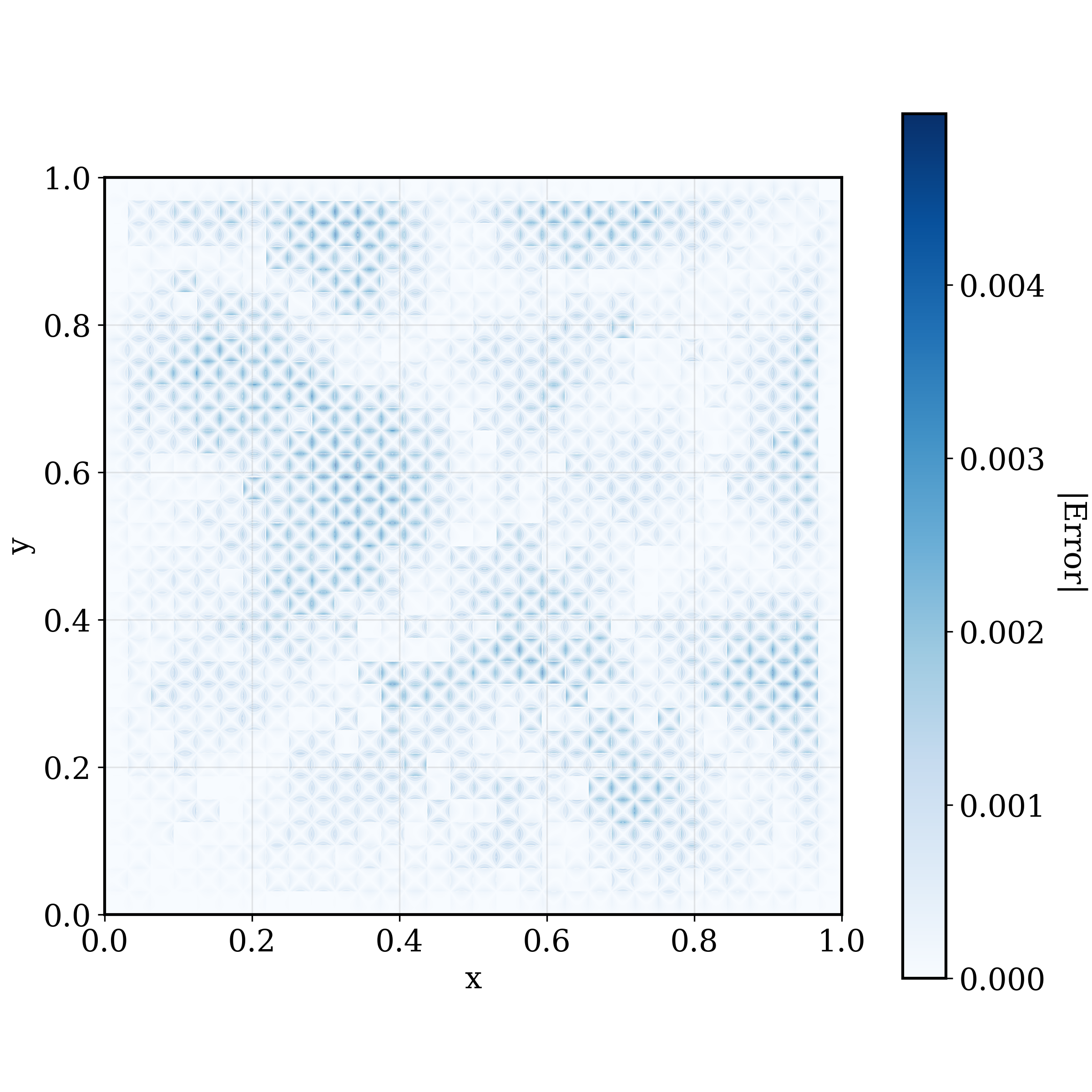}
    \caption{Error Field ($n = 1$)}
    \label{fig:ex3_hybrid_err_one}
\end{subfigure}
\hfill
\begin{subfigure}[b]{0.32\textwidth}
    \includegraphics[width=\textwidth]{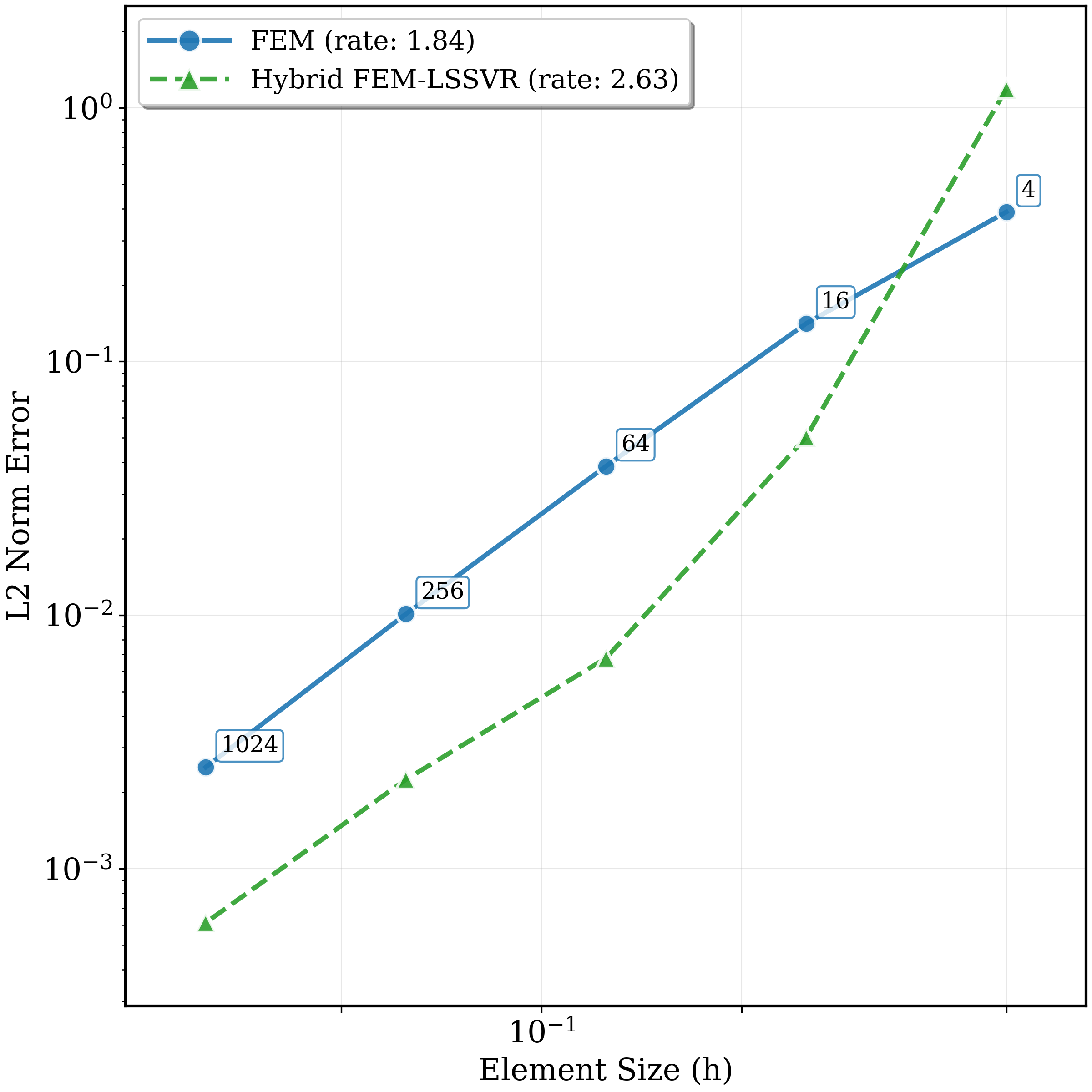}
    \caption{Convergence Analysis ($n = 1$)}
    \label{fig:ex3_convergence_one}
\end{subfigure}

\caption{Example 3 - 2D Helmholtz Equation with Dirichlet BC results showing hybrid FEM-LSSVR solution, error field, and convergence analysis for $n=\frac{1}{2}$ (top row) and $n=1$ (bottom row).}
\label{fig:ex3_results}
\end{figure}

In Figure~\ref{fig:ex3_results} we have the solutions of the Helmholtz equation with Dirichlet boundary conditions,  the error distribution, and convergence plots. Although we don't see signs of stagnation in the convergence plot, we observe consistent resolution improvement, demonstrating the hybrid method's advantage over FEM.

Table~\ref{tab5} shows the results for this example with $n = \frac{1}{2}$ and $n = 1$, using 1024 elements. In the lower-frequency case ($n = \frac{1}{2}$), all methods perform reasonably well, but the hybrid method achieves the best solution accuracy, outperforming both FEM and LSSVR. However, examining the $H^1$ norm reveals that LSSVR achieves the best accuracy, while the hybrid method and FEM show comparable performance. When the frequency increases to $n = 1$, we observe more oscillatory behavior, and the hybrid method achieves significantly better solution accuracy than both FEM and LSSVR. For the $H^1$ norm, the hybrid method maintains accuracy comparable to FEM while outperforming LSSVR.

\begin{table*}[ht!]
\centering
\caption{Performance comparison of FEM, LSSVR, and hybrid FEM+LSSVR methods for solving the 2D Helmholtz equation with Dirichlet boundary conditions. \label{tab5}}
\begin{tabular*}{\textwidth}{@{\extracolsep\fill}lccc@{\extracolsep\fill}}
\toprule
\textbf{Description/Metric} & \textbf{FEM} & \textbf{LSSVR} & \textbf{FEM+LSSVR} \\
\midrule
$n = 1/2$ \\
\midrule
Basis Order & 1 & -- & 1 \\
Kernel Order & -- & 6 & 4 \\
\midrule
Computation Time (s) & 0.067950 & 0.004089 & 5.088132 \\
Evaluation Time (s)  & 0.324757 & 0.121604 & 2.575954 \\
\midrule
Solution $\varepsilon_{\text{rel}}$ & 1.681772$\times$10$^{-3}$ & 4.503840$\times$10$^{-3}$ & 2.611473$\times$10$^{-4}$ \\
Relative H1 Error & 1.018262$\times$10$^{-2}$ & 6.696115$\times$10$^{-3}$ & 1.388744$\times$10$^{-2}$ \\
\midrule
$n = 1$ \\
\midrule
Computation Time (s) & 0.074403 & 0.010871 & 5.318305 \\
Evaluation Time (s)  & 0.455352 & 0.219404 & 2.496948 \\
\midrule
Solution $\varepsilon_{\text{rel}}$ & 6.705734$\times$10$^{-3}$ & 5.834134$\times$10$^{-2}$ & 1.773546$\times$10$^{-3}$ \\
Relative H1 Error & 2.164528$\times$10$^{-2}$ & 7.748434$\times$10$^{-2}$ & 2.816457$\times$10$^{-2}$ \\
\bottomrule
\end{tabular*}
\end{table*}

\subsection{2D Helmholtz equation with Neumann boundary conditions}
As a final example, we considered the 2D Helmholtz equation \eqref{eq:helmholtz} with Neumann boundary conditions. The equation for this specific verification test, which uses a complex wavenumber $\xi$, is given by:
\begin{equation*}
\Delta u + \xi^2 u = -f(x,y)
.
\end{equation*}
The source function is defined as
\begin{equation*}
f(x,y)=(\lambda-\xi^2)u_{\text{exact}}(x,y), \quad(x,y)\in[0,1]\times[0,1]
\end{equation*}

and homogeneous Neumann boundary conditions read
\begin{equation*}
\nabla u \cdot \nu = 0 \quad \text{on} \quad  \partial\Omega,
\end{equation*}

where $\nu$ is the outward pointing normal vector. The exact solution is
\begin{equation*}
u_{\text{exact}}(x,y)=\cos{(m \pi x)}\cos{(n \pi y)}.
\end{equation*}

Here, $\xi$ is the complex wavenumber defined by $\xi = \xi_r + i \xi_i$, where $\xi_r = \sqrt{\lambda}$ is the real part, and $\xi_i = 0.1$ is the damping factor. $\lambda$ is the eigenvalue, $\lambda = \pi^2(m^2+n^2)$, derived from the mode numbers $m$ and $n$. This problem is solved for two cases, ($m=1, n=1$) and ($m=2, n=2$). This setup allows for the evaluation of the methods' performance across two different wavelengths. For the LSSVR method, we used $5 \times 5$ uniformly distributed collocation points with $\gamma = 10^{4}$, and for the hybrid method, we used the same number in each element.

Table~\ref{tab7} presents the results for both frequencies, using $1024$ elements.
At low frequency (Mode $m=1, n=1$, corresponding to $\xi = \pi\sqrt{2} + 0.1i$), the hybrid method achieves the best solution accuracy, outperforming both FEM and LSSVR. For the $H^1$ norm, the hybrid method shows substantial improvement over both FEM and LSSVR. At high frequency (Mode $m=2, n=2$, corresponding to $\xi = 2\pi\sqrt{2} + 0.1i$), the hybrid FEM+LSSVR approach gives the best results across all metrics. The hybrid method shows significantly better solution accuracy than FEM and LSSVR, with particularly strong performance compared to LSSVR. For the $H^1$ norm, the hybrid method maintains superior accuracy compared to both methods. Compared to PINNs, where the training time increases when switching from a Dirichlet to a Neumann boundary condition \cite{schoder2025physics}, we identified a similar computational performance using the hybrid FEM+LSSVR for both types of boundary conditions. The behavior of the error metric is also recovered for the fourth test case, indicating that the hybrid method works well across these test cases.

\begin{figure}[ht!]
\centering
\begin{subfigure}[b]{0.32\textwidth}
    \centering
    \includegraphics[width=\textwidth]{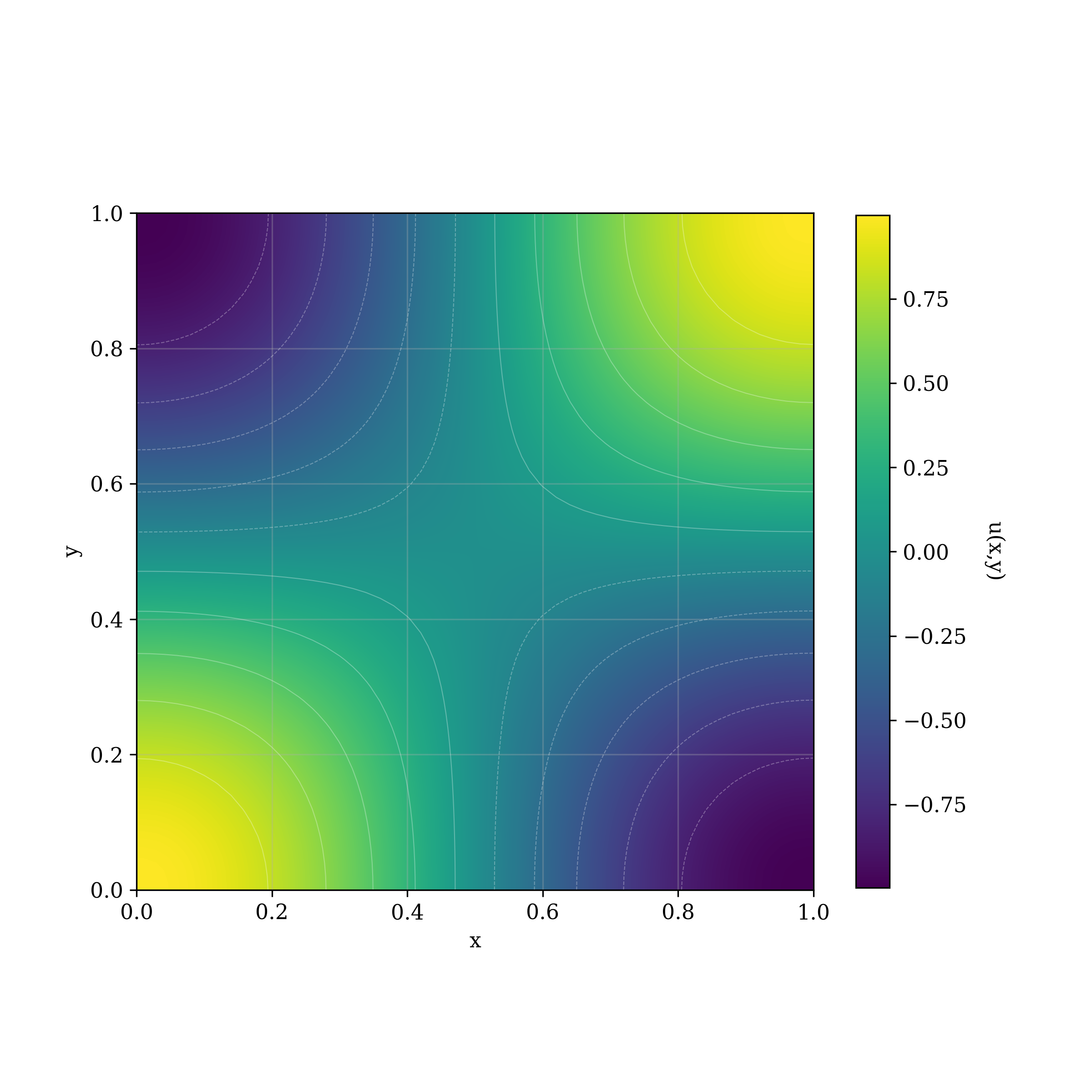}
    \caption{Hybrid Solution ($m=1, n=1$)}
    \label{fig:ex4_hybrid_sol_half}
\end{subfigure}
\hfill
\begin{subfigure}[b]{0.32\textwidth}
    \centering
    \includegraphics[width=\textwidth]{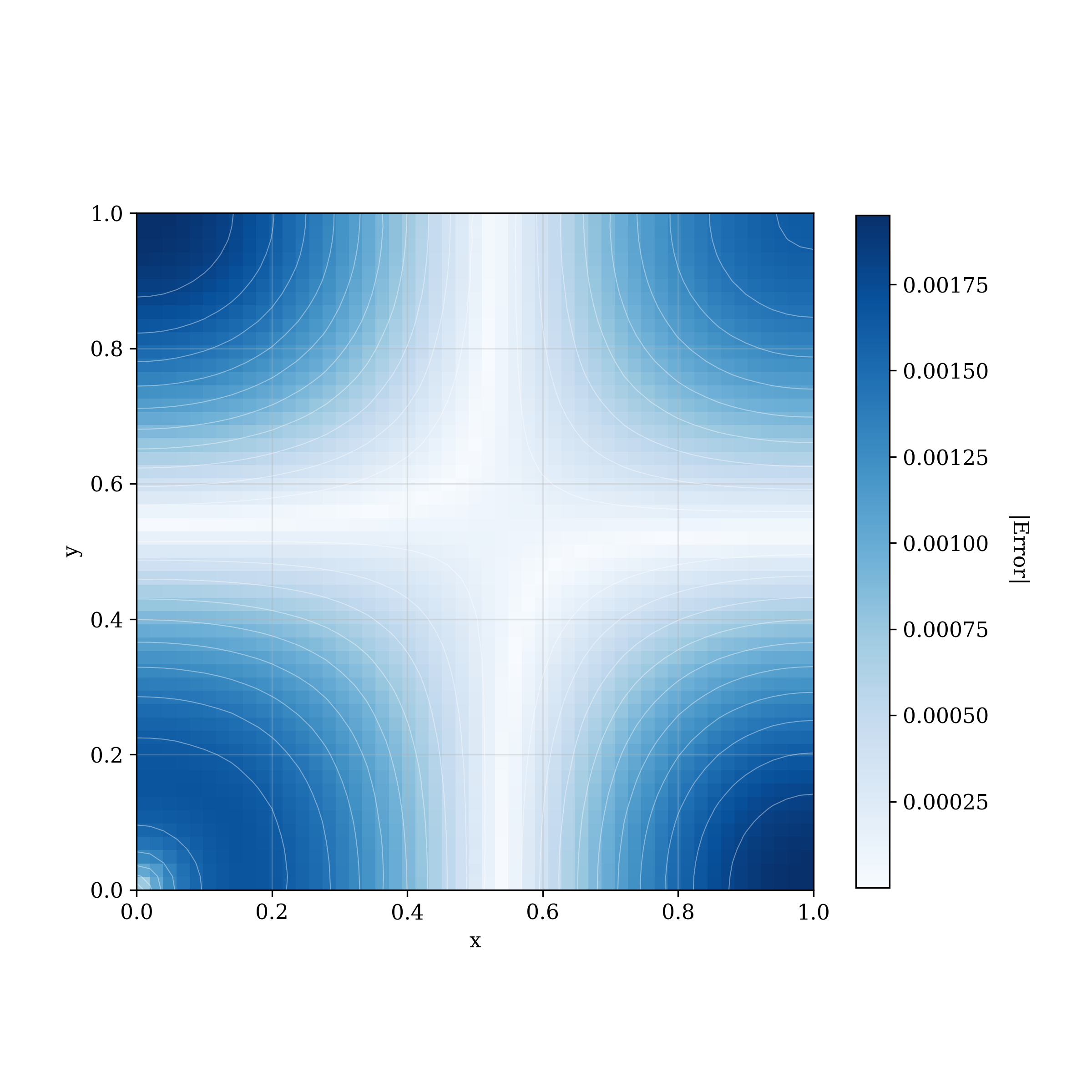}
    \caption{Error Field ($m=1, n=1$)}
    \label{fig:ex4_hybrid_err_half}
\end{subfigure}
\hfill
\begin{subfigure}[b]{0.32\textwidth}
    \centering
    \includegraphics[width=\textwidth]{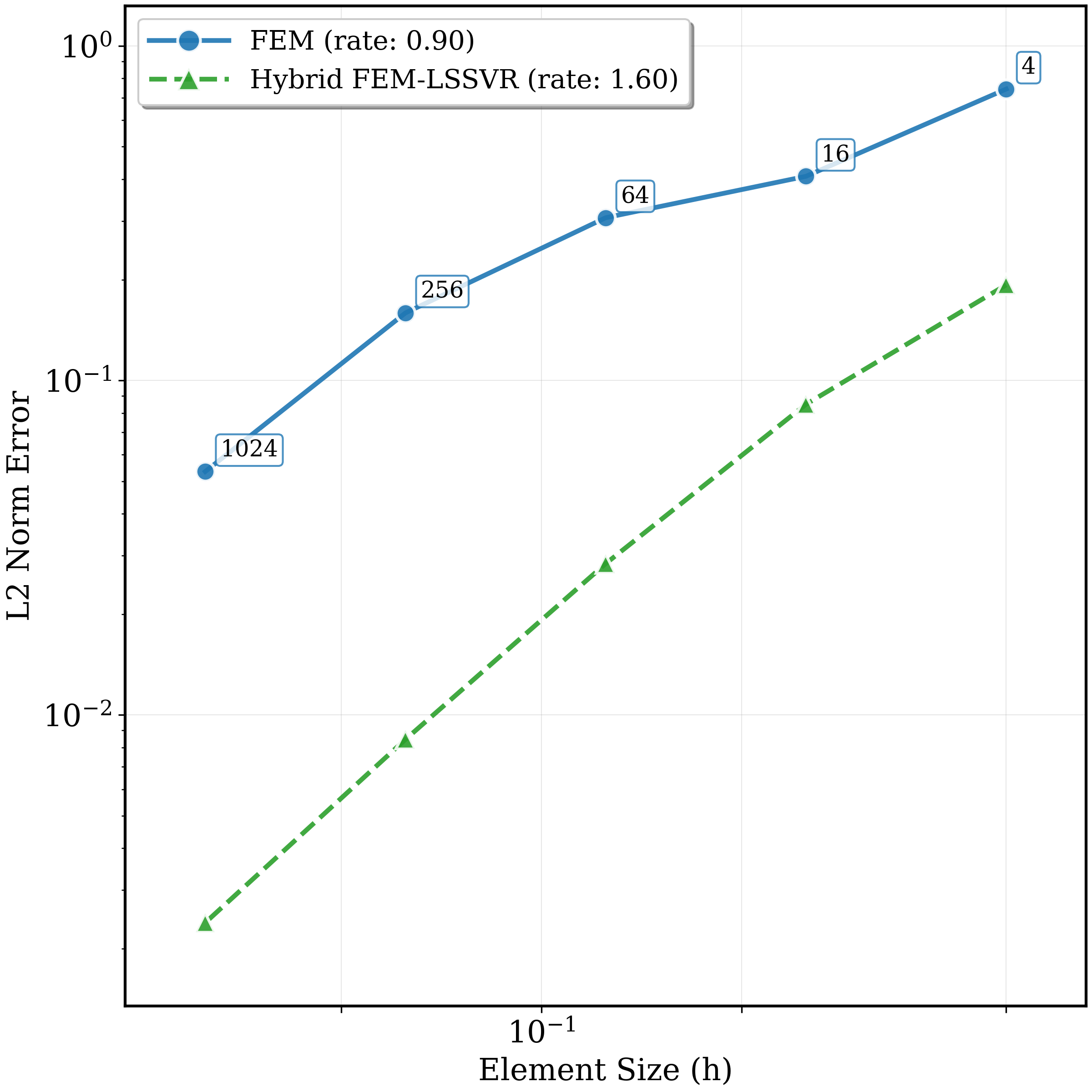}
    \caption{Convergence Analysis ($m=1, n=1$)}
    \label{fig:ex4_convergence_half}
\end{subfigure}
\vspace{0.5cm} 
\begin{subfigure}[b]{0.32\textwidth}
    \centering
    \includegraphics[width=\textwidth]{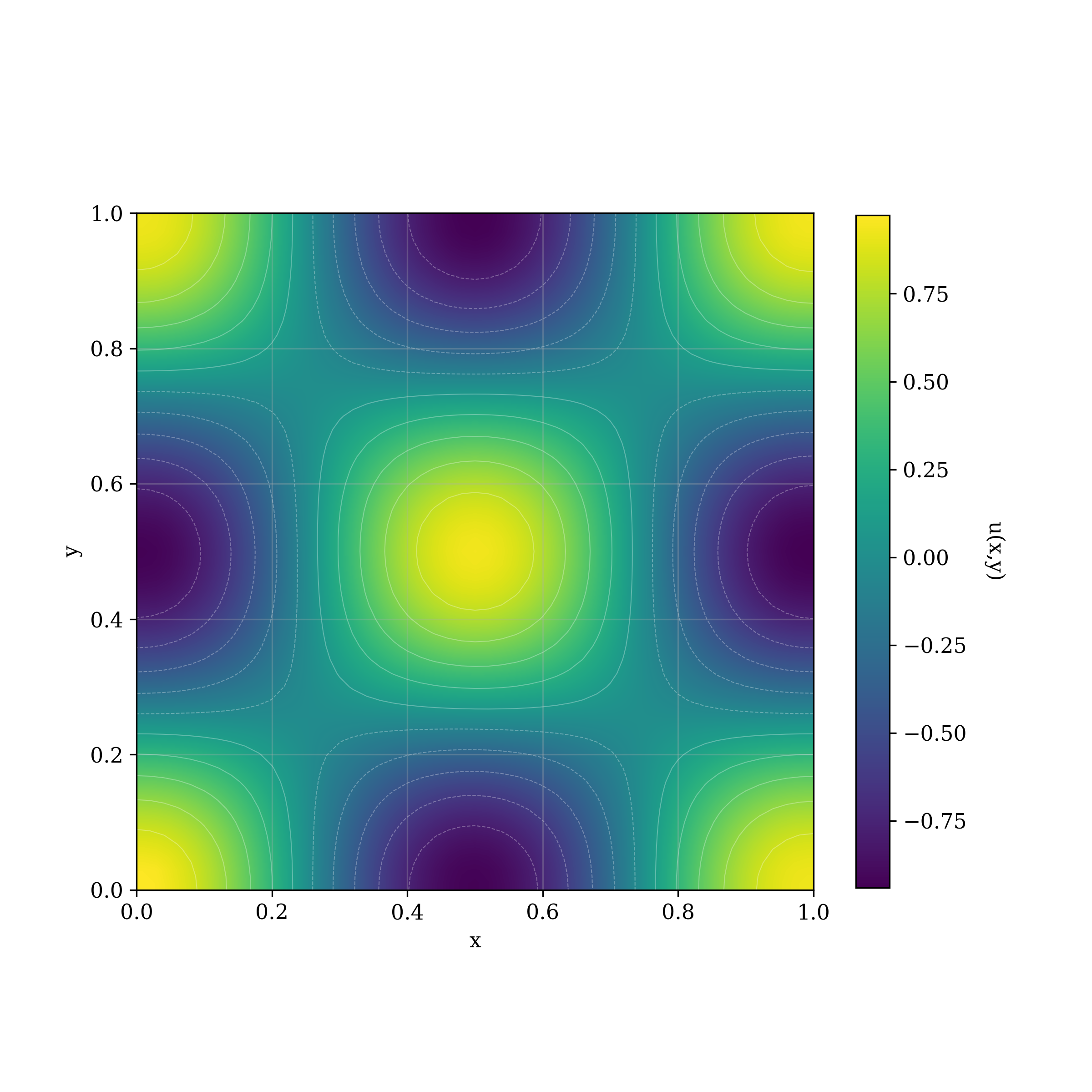}
    \caption{Hybrid Solution ($m=2, n=2$)}
    \label{fig:ex4_hybrid_sol_one}
\end{subfigure}
\hfill
\begin{subfigure}[b]{0.32\textwidth}
    \centering
    \includegraphics[width=\textwidth]{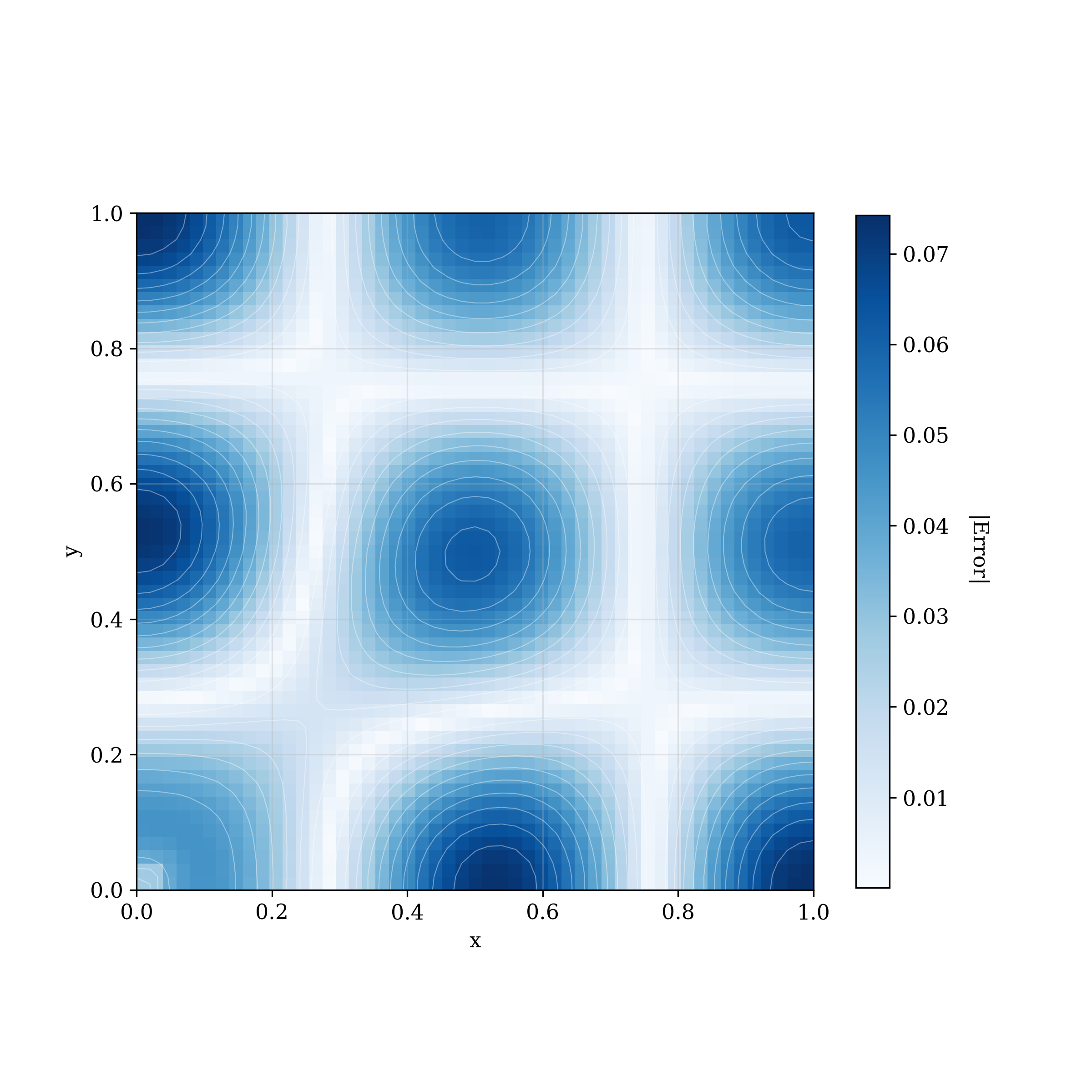}
    \caption{Error Field ($m=2, n=2$)}
    \label{fig:ex4_hybrid_err_one}
\end{subfigure}
\hfill
\begin{subfigure}[b]{0.32\textwidth}
    \centering
    \includegraphics[width=\textwidth]{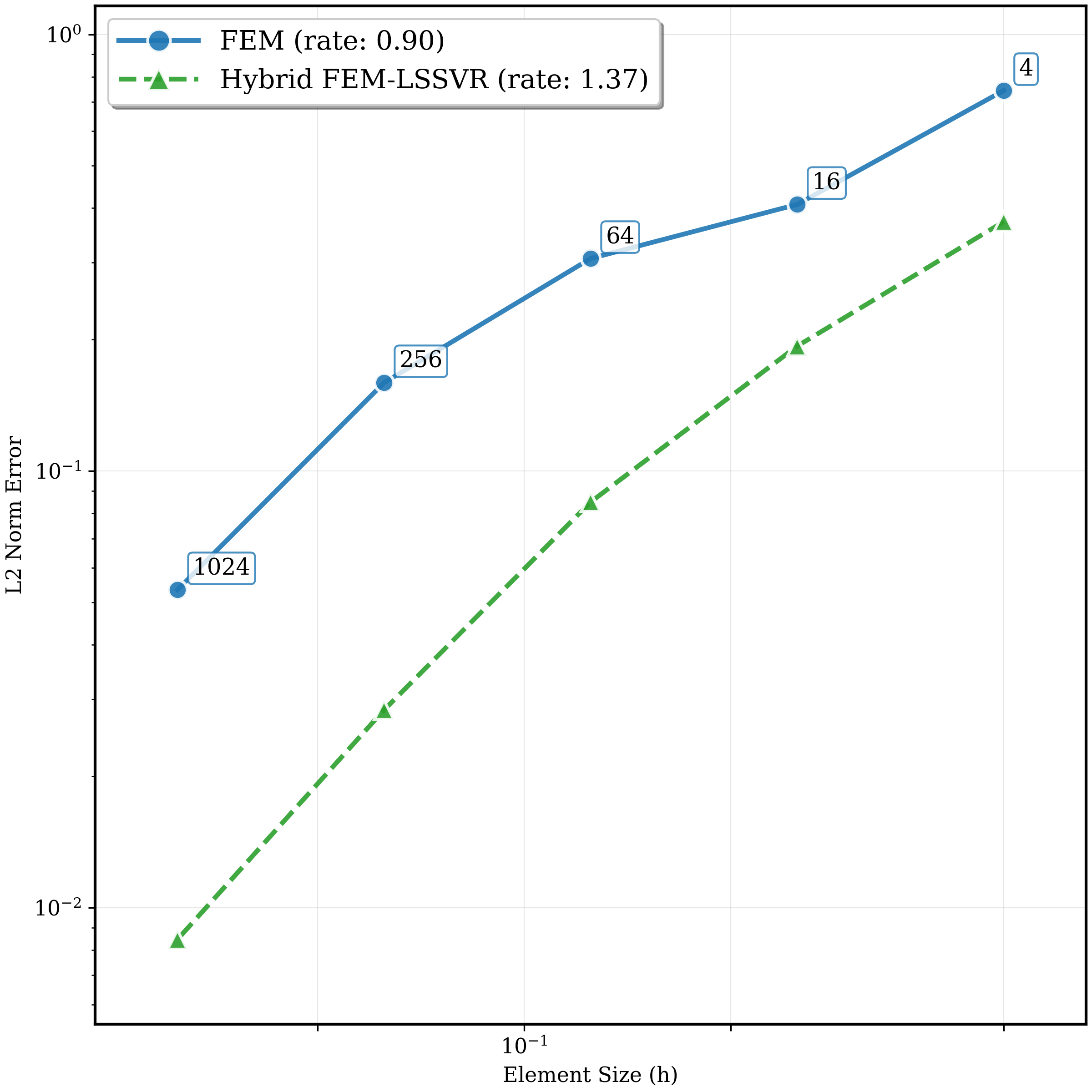}
    \caption{Convergence Analysis ($m=2, n=2$)}
    \label{fig:ex4_convergence_one}
\end{subfigure}
\caption{Example 4 - 2D Helmholtz Equation with Neumann BC results showing hybrid FEM-LSSVR solution, error field, and convergence analysis for Mode $(1, 1)$ (top row, $\xi = \pi\sqrt{2} + 0.1i$) and Mode $(2, 2)$ (bottom row, $\xi = 2\pi\sqrt{2} + 0.1i$).}
\label{fig:ex4_results}
\end{figure}

Figure~\ref{fig:ex4_results} visualizes the solution of the Helmholtz equation with Neumann boundary conditions for both frequency cases. The plots show that the hybrid method captures those cosine-based oscillatory patterns beautifully. The error plots demonstrate accurate solutions throughout the domain, even right up to the challenging Neumann boundary conditions where derivatives must be satisfied precisely.

\begin{table*}[ht!]
\centering
\caption{Performance comparison of FEM, LSSVR, and hybrid FEM+LSSVR methods for solving the 2D Helmholtz equation with Neumann boundary conditions. \label{tab7}}
\begin{tabular*}{\textwidth}{@{\extracolsep\fill}lccc@{\extracolsep\fill}}
\toprule
\textbf{Description/Metric} & \textbf{FEM} & \textbf{LSSVR} & \textbf{FEM+LSSVR} \\
\midrule

\multicolumn{4}{c}{$m = 1,\; n = 1$} \\
\midrule
Basis Order & 1 & -- & 1 \\
Kernel Order & -- & 6 & 4\\
\midrule
Computation Time (s) & $0.020483$ & $1.420839$ & $3.1046$ \\
Evaluation Time (s) & $0.237870$ & $0.045936$ & $0.8211$ \\
\midrule
Solution $\varepsilon_{\text{rel}}$ & $9.874428\times10^{-3}$ & $7.964658\times10^{-3}$ & $1.895733\times10^{-3}$ \\
Relative H1 Error & $1.586312\times10^{-2}$ & $8.059638\times10^{-3}$ & $1.849256\times10^{-3}$ \\
\midrule

\multicolumn{4}{c}{$m = 2,\; n = 2$} \\
\midrule
Computation Time (s) & $0.012223$ & $2.688048$ & $3.1509$ \\
Evaluation Time (s) & $0.227435$ & $0.309324$ & $0.5294$ \\
\midrule
Solution $\varepsilon_{\text{rel}}$ & $1.059148\times10^{-1}$ & $9.936436\times10^{-1}$ & $6.636463\times10^{-2}$ \\
Relative H1 Error & $1.096398\times10^{-1}$ & $9.936436\times10^{-1}$ & $6.521930\times10^{-2}$ \\
\bottomrule
\end{tabular*}
\end{table*}

\section{Conclusion}\label{sec:Conclusion}
In this paper, we combined the FEM and LSSVR to take advantage of both. We show that FEM provides the basic solution points, and LSSVR then improves the accuracy between those points, which is helpful when analyzing data from FEM or real-world measurements.
Numerical experiments were conducted across various equation types, domain sizes, and complexity parameters, confirming the effectiveness and reliability of the proposed method. We discussed the limitations of the method, resulting in future challenges that are addressed. Specifically, the accuracy of the solution depends on the accuracy of the provided input points posing one limiting factor for providing super-resolution when relying on numerical simulation results. Furthermore, the smoothness over the element boundaries and across the element edges will be investigated in the future.

The hybrid FEM-LSSVR method uses the geometrical definitions and boundary condition handling of the FEM, while LSSVR with Legendre polynomial kernels is used to interpolate throughout the domain. The hybrid method improves the accuracy of the solution when FEM using low-order polynomial basis functions fails to provide a sufficient accurate interpolated solution. Dirichlet and Neumann boundary conditions can be handled well, as demonstrated. The LSSVR approach adds the PDE as an equality constraint rather than including it as a term in the loss function (as it is frequently done using PINNs). In contrast to PINNs, we do not suffer from local minima issues and do not require hyperparameter tuning to balance competing objectives inside the loss function. The formulation ensures global optimality and PDE satisfaction by reformulating the problem as a convex quadratic programming problem.

The convergence analysis showed the relative error of the solution field and its derivatives to the analytic solutions. Both can be modeled accurately and converged by the hybrid formulation. 
Furthermore, the method can be used to enhance any solution obtained by low-order FEM and increase the accuracy. The solution gradients are derived using analytic differentiation with Legendre polynomials. The practical impact of this work extends beyond numerical simulation to both forward and inverse problem solving, where unknown parameters or source terms can be estimated from sparse observational data. The method can be applied to super-resolve sparse or under-resolved experimental datasets, prescribing the underlying physics by the governing PDE physics between measurement points. This transforms limited resolution of experimental data into high-resolution field representations that satisfy the underlying physical equations. 
In summary, the hybrid approach is an automated method that adds linear PDEs as conditions in the optimizer, eliminating the need for an additional hyperparameter that balances a PDE loss against a data loss.



\bmsection*{Author contributions}
M.B was drafting the manuscript and structuring/performing the formal analysis and simulations. P.R. was supporting with the FEM simulation. P.R. and M.K. were reviewing the article. S.S. was co-writing and reviewing the manuscript,  structuring the formal analysis and simulations.

\bmsection*{Acknowledgments}
Maryam Babaei, and Stefan Schoder would like to thank the Austrian Federal Ministry for Innovation, Mobility and Infrastructure and the Austrian Research Promotion Agency (FFG) for funding the bioCOMP4acoustics project under the TakeOff programme line (project number 913972). Stefan Schoder acknowledges partial funding support by the COMET project ECHODA (Energy Efficient Cooling and Heating of Domestic Appliances). ECHODA is funded within the framework of COMET - Competence Centers for Excellent Technologies by BMK, BMAW, the province of Styria as well as SFG. The COMET programme is managed by FFG. Supported by TU Graz Open Access Publishing Fund.

\bmsection*{Financial disclosure}

None reported.

\bmsection*{Conflict of interest}

The authors declare no potential conflict of interests.

\bmsection*{Data availability}

The source code and data are available at:
\url{https://github.com/maryambabaei/hybrid-FEM-LSSVR}

\bibliography{wileyNJD-AMA}

\bmsection*{Supporting information}

Additional supporting information may be found in the
online version of the article at the publisher’s website.

\end{document}